\documentclass[a4paper,12pt]{amsart}
\usepackage{mathrsfs}
\usepackage{amsfonts}
\usepackage{amssymb}
\usepackage{ifthen}
\usepackage{amscd}
\usepackage{amsxtra}
\usepackage{graphicx}
\usepackage{color}
\nonstopmode \numberwithin{equation}{section}
\newtheorem{thm}{Theorem}[section]
\newtheorem{lem}{Lemma}[section]
\newtheorem{cor}{Corollary}[section]

\newtheorem{step}{Step}[section]

\newtheorem{cl}{Claim}[section]
\newtheorem{ca}{Case}
\newtheorem{sca}{Subcase}[section]
\newtheorem{scl}{Subclaim}[section]
\newtheorem{conj}{Conjecture}

\theoremstyle{definition}
\newtheorem{defn}{Definition}[section]

\newtheorem{op}[equation]{Open Problem}
\newtheorem{ques}[equation]{Question}
\newtheorem{rem}{Remark}[section]
\newtheorem{exam}[equation]{Example}

\newcounter {own}
\def\theown {\thesection       .\arabic{own}}

\newenvironment{pf}[1][]{%
 \vskip 3mm
 \noindent
 \ifthenelse{\equal{#1}{}}%
  {{\slshape Proof. }}%
  {{\slshape #1.} }%
 }%
{\qed\bigskip}

\newcounter{alphabet}

\newenvironment{Thm}[1][]{\refstepcounter{alphabet}%
\bigskip%
\noindent%
{\bf Theorem \Alph{alphabet}}%
\ifthenelse{\equal{#1}{}}{}{ (#1)}%
{\bf .} \itshape}{\vskip 8pt}


\newenvironment{Lem}[1][]{\refstepcounter{alphabet}%
\bigskip%
\noindent%
{\bf Lemma \Alph{alphabet}}%
{\bf .} \itshape}{\vskip 8pt}

\newcommand{\IR}{{\mathbb R}}

\newcommand{\Aut}{{\operatorname{Aut}}}

\def\be{\begin{equation}}
\def\ee{\end{equation}}

\newcommand{\bee}{\begin{enumerate}}
\newcommand{\eee}{\end{enumerate}}

\newcommand{\blem}{\begin{lem}}
\newcommand{\elem}{\end{lem}}
\newcommand{\bthm}{\begin{thm}}
\newcommand{\ethm}{\end{thm}}
\newcommand{\bcor}{\begin{cor}}
\newcommand{\ecor}{\end{cor}}
\newcommand{\beg}{\begin{exam}}
\newcommand{\eeg}{\end{exam}}
\newcommand{\begs}{\begin{examples}}
\newcommand{\eegs}{\end{examples}}
\newcommand{\bdefe}{\begin{defn}}
\newcommand{\edefe}{\end{defn}}
\newcommand{\bprob}{\begin{pro}}
\newcommand{\eprob}{\end{pro}}
\newcommand{\bques}{\begin{ques}}
\newcommand{\eques}{\end{ques}}
\newcommand{\bei}{\begin{itemize}}
\newcommand{\eei}{\end{itemize}}
\newcommand{\bcon}{\begin{conj}}
\newcommand{\econ}{\end{conj}}
\newcommand{\bop}{\begin{op}}
\newcommand{\eop}{\end{op}}

\newcommand{\bca}{\begin{ca}}
\newcommand{\eca}{\end{ca}}
\newcommand{\bsca}{\begin{sca}}
\newcommand{\esca}{\end{sca}}

\newcommand{\bcl}{\begin{cl}}
\newcommand{\ecl}{\end{cl}}

\newcommand{\bst}{\begin{step}}
\newcommand{\est}{\end{step}}

\newcommand{\bscl}{\begin{scl}}
\newcommand{\escl}{\end{scl}}

\newcommand{\bcons}{\begin{conjs}}
\newcommand{\econs}{\end{conjs}}
\newcommand{\bprop}{\begin{propo}}
\newcommand{\eprop}{\end{propo}}
\newcommand{\br}{\begin{rem}}
\newcommand{\er}{\end{rem}}
\newcommand{\brs}{\begin{rems}}
\newcommand{\ers}{\end{rems}}
\newcommand{\bo}{\begin{obser}}
\newcommand{\eo}{\end{obser}}
\newcommand{\bos}{\begin{obsers}}
\newcommand{\eos}{\end{obsers}}
\newcommand{\bpf}{\begin{pf}}
\newcommand{\epf}{\end{pf}}
\newcommand{\ba}{\begin{array}}
\newcommand{\ea}{\end{array}}
\newcommand{\beq}{\begin{eqnarray}}
\newcommand{\beqq}{\begin{eqnarray*}}
\newcommand{\eeq}{\end{eqnarray}}
\newcommand{\eeqq}{\end{eqnarray*}}

\newcommand{\ds}{\displaystyle}

\newcounter{minutes}
\divide\time by 60
\newcounter{hours}
\multiply\time by 60 \addtocounter{minutes}{-\time}

\begin{document}
\bibliographystyle{amsplain}
\title [] {The Heinz type inequality, Bloch type theorem and Lipschitz characteristic of polyharmonic mappings}

\def\thefootnote{}
\footnotetext{ \texttt{\tiny File:~\jobname .tex,
          printed: \number\day-\number\month-\number\year,
          \thehours.\ifnum\theminutes<10{0}\fi\theminutes}
} \makeatletter\def\thefootnote{\@arabic\c@footnote}\makeatother

\author{Shaolin Chen}
\address{S. L. Chen,    College of Mathematics and
Statistics, Hengyang Normal University, Hengyang, Hunan 421008,
People's Republic of China} \email{mathechen1982@hynu.edu.cn}




\begin{abstract}

Suppose that $f$  satisfies the following: $(1)$ the polyharmonic
  equation $\Delta^{m}f=\Delta(\Delta^{m-1}
f)$$=\varphi_{m}$ $(\varphi_{m}\in
\mathcal{C}(\overline{\mathbb{B}^{n}},\mathbb{R}^{n}))$, (2) the
boundary conditions
$\Delta^{0}f=\varphi_{0},\Delta^{1}f=\varphi_{1},~\ldots,~\Delta^{m-1}f=\varphi_{m-1}$
on $\mathbb{S}^{n-1}$ ($\varphi_{j}\in
\mathcal{C}(\mathbb{S}^{n-1},\mathbb{R}^{n})$ for
$j\in\{0,1,\ldots,m-1\}$ and $\mathbb{S}^{n-1}$ denotes the boundary
of the unit ball $\mathbb{B}^{n}$), and $(3)$ $f(0)=0$, where
$n\geq3$ and $m\geq1$ are integers.  Initially, we prove a Schwarz type lemma and use  it to obtain a Heinz type
inequality of mappings satisfying the
polyharmonic equation with the above Dirichlet boundary value
conditions. Furthermore, we establish a Bloch type theorem of mappings satisfying
 the above
polyharmonic equation, which gives an
answer to an open problem in \cite{CP-Hi}. Additionally,
 we  show that if $f$ is a $K$-quasiconformal self-mapping of  $\mathbb{B}^{n}$ satisfying
 the above
polyharmonic equation, then $f$ is Lipschitz continuous, and the
Lipschitz constant is asymptotically sharp as $K\to 1^{+}$ and
$\|\varphi_{j}\|_{\infty}\to 0^{+}$ for $j\in\{1,\ldots,m\}$.
\end{abstract}




\maketitle \pagestyle{myheadings} \markboth{S.  L. Chen} { The Heinz type inequality, Bloch type theorem and Lipschitz characteristic}

\section{Preliminaries and statements of main results}\label{csw-sec1}
For an  integer $n\geq2$, let $\mathbb{R}$ and $\mathbb{R}^{n}$ be
the set of real numbers
 and the usual real vector space of dimension $n$, respectively. Sometimes it is convenient to
identify each point $x\in\mathbb{R}^{n}$ with
an $n\times 1$ column matrix so that $x=(x_{1},\ldots,x_{n})'$,
where $'$ denotes the transposition  of a matrix.
For $y=(y_{1},\ldots,y_{n})'$ and $x\in\mathbb{R}^{n}$, we define
the Euclidean inner product $\langle \cdot ,\cdot \rangle$ by
$\langle x,y\rangle=\sum_{k=1}^{n}x_{k}y_{k}$ so that the Euclidean
length of $x$ is defined by
$$|x|=\langle x,x\rangle^{1/2}=\left(\sum_{k=1}^{n}|x_{k}|^{2}\right)^{1/2}.
$$
Denote a ball in $\mathbb{R}^{n}$ with center
$x_{0}\in\mathbb{R}^{n}$ and radius $r$ by
$\mathbb{B}^{n}(x_{0},r)$. In particular, let
$\mathbb{B}^{n}:=\mathbb{B}^{n}(0,1)$ and
$\mathbb{S}^{n-1}:=\partial\mathbb{B}^{n}$.
 For
$n_{1}\in\mathbb{N}:=\{1,2,\ldots\}$ and
$k\in\mathbb{N}_0=\mathbb{N}\cup\{0\}$, we denote by
$\mathcal{C}^{k}(\Omega_{1},\Omega_{2})$ the set of all $k$-times
continuously differentiable functions from $\Omega_{1}$ into
$\Omega_{2}$, where $\Omega_{1}$ and $\Omega_{2}$ are subsets of
$\mathbb{R}^{n}$ and $\mathbb{R}^{n_{1}}$, respectively. In
particular, let
$\mathcal{C}(\Omega_{1},\Omega_{2}):=\mathcal{C}^{0}(\Omega_{1},\Omega_{2})$,
the set of all continuous functions of $\Omega_{1}$ into
$\Omega_{2}$. For
$f=(f_{1},\ldots,f_{n_{1}})'\in\mathcal{C}^{1}(\Omega_{1},\Omega_{2})$,
we denote the derivative $D_{f}$ of $f$ by

$$D_{f}=\left(\begin{array}{cccc}
\ds D_{1}f_{1}\; \cdots\;
 D_{n}f_{1}\\[4mm]
\vdots\;\; \;\;\cdots\;\;\;\;\vdots \\[2mm]
 \ds
D_{1}f_{n_{1}}\; \cdots\;
 D_{n}f_{n_{1}}
\end{array}\right), \quad D_{j}f_{i}(x)=\frac{\partial f_{i}(x)}{\partial x_j}.
$$

In particular, if $n=n_{1}$, the Jacobian of $f$ is defined by
$J_{f}=\det D_{f}$ and the Laplacian of
$f\in\mathcal{C}^{2}(\Omega_{1},\Omega_{2})$ is defined by
$$\Delta f=\sum_{k=1}^{n}D_{kk}f.$$

For an $n_{1}\times n$ matrix $A=(a_{ij})_{n_{1}\times n}$ with
$n_{1}\in\mathbb{N}$ and $n\geq2$, the operator norm of $A$ is
defined by
$$|A|=\sup_{x\in\mathbb{R}^{n},x\neq 0}\frac{|Ax|}{|x|}=\max\{|A\theta|:\,
\theta\in\mathbb{S}^{n-1}\},
$$
and the matrix function $l(A)$ is defined by
$$l(A)=\inf\{|A\theta|:~\theta\in\mathbb{S}^{n-1}\}.$$


\subsection{Polyharmonic equation}
For $n\geq3$ and $x,y\in\mathbb{R}^{n}\backslash\{0\}$, we define
$x^{\ast}=x/|x|$, $y^{\ast}=y/|y|$ and let
$$[x,y]:=\left|y|x|-x^{\ast}\right|=\left|x|y|-y^{\ast}\right|.
$$
Also, for $x,y\in\mathbb{B}^{n}$ with $x\neq y$,
we use $G(x,y)$ to denote the {\it Green function}:
\be\label{eq-ex0}
G(x,y)=c_{n}\left(\frac{1}{|x-y|^{n-2}}-\frac{1}{[x,y]^{n-2}}\right),
\ee where $c_{n}=1/[(n-2)\omega_{n-1}]$ and
$\omega_{n-1}=2\pi^{\frac{n}{2}}/\Gamma\big(\frac{n}{2}\big)$
denotes the area of $\mathbb{S}^{n-1}$. The {\it
Poisson kernel} $P:\,\mathbb{B}^{n}\times
\mathbb{S}^{n-1}\rightarrow \IR$ is defined by
$$P(x,\zeta)=\frac{1-|x|^{2}}{|x-\zeta|^{n}}.
$$
We use
$$\nabla =\left( \frac{\partial }{\partial x_1}, \ldots, \frac{\partial }{\partial x_n} \right)
$$ to denote the {\it gradient}.

Of particular interest for our investigation is the following {\it
polyharmonic} equation:

\be\label{eq-ch-1} \Delta^{m}f=\Delta(\Delta^{m-1}
f)=\varphi_{m}~\mbox{in}~\mathbb{B}^{n}\ee with the following
associated {\it Dirichlet boundary value condition}:

\be\label{eq-ch-2}
\Delta^{0}f=\varphi_{0},~\Delta^{1}f=\varphi_{1},~\ldots,~\Delta^{m-1}f=\varphi_{m-1}~\mbox{on}~\mathbb{S}^{n-1},\ee
where $m,~n_{1}\in\mathbb{N}$,  $\Delta^{0}f:=f$, $\Delta^{1}
f:=\Delta f$,
$\varphi_{m}\in\mathcal{C}(\mathbb{B}^{n},\mathbb{R}^{n_{1}})$, and
$\varphi_{k}\in \mathcal{C}(\mathbb{S}^{n-1},\mathbb{R}^{n_{1}})$
for $k\in\{0,1,\ldots,m-1\}$.  Here the boundary condition in
(\ref{eq-ch-2}) are interpreted in the following distributional
sense. For some fixed $r\in(0,1)$, let
$f_{r}(x)=f(rx),~x\in\mathbb{B}^{n}$. Then for any
$\zeta\in\mathbb{S}^{n-1}$,
$$\lim_{r\rightarrow1^{-}}\Delta^{j-1}f_{r}(\zeta)=\varphi_{j-1}(\zeta),$$
where $j\in\{1,\ldots,m\}$.

(I) If $m=1$, then all solutions to the equation (\ref{eq-ch-1}) satisfying
(\ref{eq-ch-2}) are given by
\be\label{po-1}f(x)=P[\varphi_{0}](x)-G_{1}[\varphi_{1}](x),~x\in\mathbb{B}^{n},\ee
where
$$P[\varphi_{0}](x)=\int_{\mathbb{S}^{n-1}}P(x,\zeta)\varphi_{0}(\zeta)d\sigma(\zeta)$$
and \be\label{po-2}G_{1}[\varphi_{1}](x)=\int_{\mathbb{B}^{n}}G(x,y_{1})\varphi_{1}(y_{1})dV(y_{1}).\ee
Here $d\sigma$ denotes the normalized Lebesgue
surface measure on $\mathbb{S}^{n-1}$ and $dV$ is the Lebesgue
volume measure on $\mathbb{B}^{n}$.

(II) If $m\geq2$, then, by \cite[p.~118--120]{Ho}  and the iterative procedure, we see that
all solutions to the equation (\ref{eq-ch-1}) satisfying
(\ref{eq-ch-2}) are given by

\be\label{eq-ch-3.0}f(x)=P[\varphi_{0}](x)+\sum_{k=1}^{m}(-1)^{k}G_{k}[\varphi_{k}](x),~x\in\mathbb{B}^{n},\ee
where

\beq\label{eq-ch-3.1} G_{k}[\varphi_{k}](x)&=&
\int_{\mathbb{B}^{n}}\cdots\int_{\mathbb{B}^{n}}G(x,y_{1})\cdots
G(y_{k-1},y_{k})\\ \nonumber
&&\times\left(\int_{\mathbb{S}^{n-1}}P(y_{k},\zeta)\varphi_{k}(\zeta)d\sigma(\zeta)\right)dV(y_{k})\cdots
dV(y_{1})\eeq for  $k\in\{1,\ldots,m-1\}$, and

\beq\label{eq-ch-3.2} G_{m}[\varphi_{m}](x)&=&
\int_{\mathbb{B}^{n}}\cdots\int_{\mathbb{B}^{n}}G(x,y_{1})\cdots G(y_{m-2},y_{m-1})\\
\nonumber
&&\times\left(\int_{\mathbb{B}^{n}}G(y_{m-1},y_{m})\varphi_{m}(y_{m})dV(y_{m})\right)dV(y_{m-1})\cdots
dV(y_{1}).\eeq
Moreover,
we call $f$ a polyharmonic mapping if $f$ satisfies (\ref{eq-ch-3.0}).

We refer the reader to
\cite{CCGS,GGS,May-M} etc for more discussions in this line. In
particular, if $m=1$ ($m=2$ resp.), then (\ref{eq-ch-1}) is called
the {\it Poisson equation} ( {\it biharmonic equation} resp.) (cf.
\cite{K1,K3,K-2011,Lai}).


\subsection{Main results}
Heinz in his classical paper \cite{He} showed that the following
result which is called the  Schwarz  Lemma  of harmonic mappings: If $f$ is a harmonic mapping of
the unit disk $\mathbb{D}:=\mathbb{B}^{2}$ into $\mathbb{D}$ with
$f(0)=0$, then
$$|f(z)|\leq\frac{4}{\pi}\arctan|z|.$$ Later,  Pavlovi\'c \cite[Theorem 3.6.1]{Pav1} removed
the assumption $f(0)=0$ and obtained the following sharp form

\be\label{eq-pav1}\left|f(z)-\frac{1-|z|^{2}}{1+|z|^{2}}f(0)\right|\leq\frac{4}{\pi}\arctan
|z|,\ee where $f$ is a harmonic mapping from $\mathbb{D}$ into
itself. The inequality (\ref{eq-pav1}) has been  proved
independently by Hethcote in
 \cite{Het}.
 For $n\geq3$, the
classical Schwarz lemma of harmonic mappings in $\mathbb{B}^{n}$
infers that if $f$ is a harmonic mapping of $\mathbb{B}^{n}$ into
itself satisfying $f(0)=0,$ then
$$|f(x)|\leq U(rN),
$$
where $r=|x|$, $N=(0,\ldots,0,1)'$ and $U$ is a harmonic function of
$\mathbb{B}^{n}$ into $[-1,1]$ defined by
$$U(x)=P[{X}_{S^{+}}-\mathcal{X}_{S^{-}}](x).
$$
Here $\mathcal{X}$ is the indicator function,
$S^{+}=\{x=(x_{1},\ldots,x_{n})'\in\mathbb{S}^{n-1}:~x_{n}\geq0\}$
and
$S^{-}=\{x=(x_{1},\ldots,x_{n})'\in\mathbb{S}^{n-1}:~x_{n}\leq0\}$
(see \cite{ABR}). In \cite{K5}, Kalaj showed that the following
result for harmonic mappings $f$ of $\mathbb{B}^{n}$ into itself:
\be\label{eq-K}
\left|f(x)-\frac{1-|x|^{2}}{(1+|x|^{2})^{\frac{n}{2}}}f(0)\right|\leq
U(|x|N). \ee 

The first aim of the paper is to extend  \eqref{eq-K} to mappings
satisfying the polyharmonic equation. More precisely, we shall prove
the following.

\begin{thm}\label{thm-1}
Let $n\geq3$, $m\geq2,~n_{1}\in\mathbb{N}$,
$\varphi_{m}\in\mathcal{C}(\overline{\mathbb{B}^{n}},\mathbb{R}^{n_{1}})$
and $\varphi_{k}\in
\mathcal{C}(\mathbb{S}^{n-1},\mathbb{R}^{n_{1}})$ for
$k\in\{1,\ldots,m-1\}$. If
$f\in\mathcal{C}^{2m}(\mathbb{B}^{n},\mathbb{R}^{n_{1}})\cap\mathcal{C}(\overline{\mathbb{B}^{n}},\mathbb{R}^{n_{1}})$
satisfies \eqref{eq-ch-1} with the boundary condition:
$\Delta^{1}f=\varphi_{1},~\ldots,~\Delta^{m-1}f=\varphi_{m-1}$ on
$\mathbb{S}^{n-1}$, then for $x\in\overline{\mathbb{B}^{n}}$,

\beq\label{eq-thm1}
\left|f(x)-\frac{1-|x|^{2}}{(1+|x|^{2})^{\frac{n}{2}}}P[\varphi_{0}](0)\right|
&\leq&\|\varphi_{0}\|_{\infty}U(|x|N)\\ \nonumber
&&+\sum_{k=1}^{m}\frac{\|\varphi_{k}\|_{\infty}}{2n}\left[\frac{n+4}{4n(n+2)}\right]^{k-1}(1-|x|^{2}),\eeq
where $\varphi_{0}=f|_{\mathbb{S}^{n-1}}$,
$\|\varphi_{m}\|_{\infty}=\sup_{x\in\mathbb{B}^{n}}|\varphi_{m}(x)|~
\mbox{and}~
\|\varphi_{k}\|_{\infty}=\sup_{\zeta\in\mathbb{S}^{n-1}}|\varphi_{k}(\zeta)|$
for  $k\in\{0,1,\ldots,m-1\}$.

In particular, if we choose
$f(x)=(M(1-|x|^{2(m-1)}),0,\ldots,0)\in\mathcal{C}^{2m}(\mathbb{B}^{n},\mathbb{R}^{n_{1}})\cap\mathcal{C}(\overline{\mathbb{B}^{n}},\mathbb{R}^{n_{1}})$
for $x\in\overline{\mathbb{B}^{n}}$, then
 the inequality
{\rm(\ref{eq-thm1})} is sharp in $\mathbb{S}^{n-1}$, where $M>0$ is
a constant.
\end{thm}

Let $f$ be a harmonic homeomorphism  of $\mathbb{D}$ onto itself
with $f(0)=0.$  Heinz  \cite[Ineq. (18)]{He} proved that, for any
$\theta\in[0,2\pi]$,
$$\liminf_{r\rightarrow1^{-}}|D_{f}(re^{i\theta})|\geq\frac{1}{\pi}.$$
We refer to \cite{K5} for the extensive discussion on  Heinz type
inequalities for harmonic mappings in $\mathbb{R}^{n}~(n\geq3)$. On
the applications of the Heinz type inequalities, see
\cite{Du,K-2015}. In the following, by using Theorem \ref{thm-1}, we
establish a Heinz type inequality for mappings satisfying the
polyharmonic equation.

\begin{thm}\label{thm-2}
For $n\geq3$, $m\geq2,~n_{1}\in\mathbb{N}$ and
$k\in\{1,\ldots,m-1\}$, suppose that
$\varphi_{m}\in\mathcal{C}(\overline{\mathbb{B}^{n}},\mathbb{R}^{n_{1}})$
and $\varphi_{k}\in
\mathcal{C}(\mathbb{S}^{n-1},\mathbb{R}^{n_{1}})$
 such that
$$\frac{n!\big[1+n-(n-2)F\big(\frac{1}{2},1;\frac{n+3}{2};-1\big)\big]}{2^{\frac{3n}{2}}\Gamma(\frac{n+1}{2})\Gamma(\frac{n+3}{2})}>
\sum_{k=1}^{m}\frac{\|\varphi_{k}\|_{\infty}}{n}\left[\frac{n+4}{4n(n+2)}\right]^{k-1}\left(1+\frac{1}{2^{\frac{n}{2}}}\right).$$
Let
$f\in\mathcal{C}^{2m}(\mathbb{B}^{n},\mathbb{R}^{n_{1}})\cap\mathcal{C}(\overline{\mathbb{B}^{n}},\mathbb{R}^{n_{1}})$
satisfying \eqref{eq-ch-1} with the Dirichlet boundary value
condition:
$\Delta^{1}f=\varphi_{1},~\ldots,~\Delta^{m-1}f=\varphi_{m-1}$ on
$\mathbb{S}^{n-1}$. If $f(0)=0$ and
$\lim_{r\rightarrow1^{-}}|f(r\zeta)|=1$ for some
$\zeta\in\mathbb{S}^{n-1}$, then
\beq\label{chk}
\liminf_{r\rightarrow1^{-}}\frac{|f(\zeta)-f(r\zeta)|}{1-r}&\geq&\frac{n!\big[1+n-(n-2)F\big(\frac{1}{2},1;\frac{n+3}{2};-1\big)\big]}{2^{\frac{3n}{2}}\Gamma(\frac{n+1}{2})\Gamma(\frac{n+3}{2})}\\
\nonumber
&&-\sum_{k=1}^{m}\frac{\|\varphi_{k}\|_{\infty}}{n}\left[\frac{n+4}{4n(n+2)}\right]^{k-1}\left(1+\frac{1}{2^{\frac{n}{2}}}\right),
\eeq where $\Gamma$ and $F(\cdot,\cdot;\cdot;\cdot)$ are the Gamma
function and the hypergeometric function, respectively $($see the
Section {\rm \ref{sbcsw-sec2.2}}$)$.
 In particular, if $\|\varphi_{k}\|_{\infty}=0$ for
$k\in\{1,\ldots,m\}$, then this estimate {\rm (\ref{chk})} is sharp.
\end{thm}

Next, we  discuss an issue that  is related to a classical result in
geometric function theory: the theorem of Bloch. Recall that Bloch's
theorem says that  an analytic function $f$ on the unit disk with
$|f'(0)|=1$ univalently covers a disk of radius $R$, where $R$ is a
universal constant (see \cite{Raj,Raj-2004}). However, for general
class of functions, there is no Bloch's Theorem. For example,
consider
 $f_{k}(x)=(kx_{1},x_{2}/k,
x_{3},\ldots,x_{n})'$ for $k=\{1,2,\ldots\}$, where $n\geq3$ and
$x=(x_{1},\ldots,x_{n})'\in\mathbb{B}^{n}$. It is easy to see that
each $f_{k}$ is univalent and $|f_{k}(0)|=J_{f_{k}}(0)-1=0$.
Furthermore, each  $f_{k}(\mathbb{B}^{n})$ contains no ball with
radius bigger than $1/k$. Hence, there does not exist an absolute
constant $R_{0}>0$ which can work for all $k\in\{1,2,\ldots\}$, such
that $\mathbb{B}(0,R_{0}) $ is contained in the range
$f_{k}(\mathbb{B}^{n})$. To establish analogs of the Bloch's theorem
for more general classes of functions, it is necessary to restrict
our focus on certain subclasses (see
\cite{AJK,BE,CG,CLW-2018,CP-Hi,Er,Raj,Raj-1,W}). In our next result, we
establish a Bloch type theorem for mappings satisfying the
polyharmonic equation, which gives an answer to the open problem in
\cite[Remark 1.2]{CP-Hi}.

\begin{defn}
Let $M$ be a positive constant and $n\geq3$.
\begin{enumerate}
\item[{\rm (1)}] If $m=1$ in (\ref{eq-ch-1}), then, for a given function $\varphi_{1}\in\mathcal{C}(\mathbb{B}^{n},\mathbb{R}^{n})$,  we use
$\mathcal{F}_{\varphi_{1}}^{M}$ to denote the set of all mappings
$f\in\mathcal{C}^{2}(\mathbb{B}^{n},\mathbb{R}^{n})\cap\mathcal{C}(\overline{\mathbb{B}^{n}},\mathbb{R}^{n})$
satisfying $|f|\leq M$, $|f(0)|=J_{f}(0)-1=0$ and $ \Delta
f=\varphi_{1}~\mbox{in}~\mathbb{B}^{n}$.

\item[{\rm (2)}] If $m\geq2$ in (\ref{eq-ch-1}), then, for   given functions $\varphi_{m}\in\mathcal{C}(\mathbb{B}^{n},\mathbb{R}^{n})$ and $\varphi_{1},\ldots,\varphi_{m-1}\in
\mathcal{C}(\mathbb{S}^{n-1},\mathbb{R}^{n})$, we denote by
$\mathcal{F}_{\varphi_{1},\ldots,\varphi_{m}}^{M}$ the set of all
mappings
$f\in\mathcal{C}^{2m}(\mathbb{B}^{n},\mathbb{R}^{n})\cap\mathcal{C}(\overline{\mathbb{B}^{n}},\mathbb{R}^{n})$
satisfying $|f|\leq M$, $|f(0)|=J_{f}(0)-1=0$, and \eqref{eq-ch-1}
with the following
  Dirichlet boundary value condition:

$$
\Delta^{1}f=\varphi_{1},~\ldots,~\Delta^{m-1}f=\varphi_{m-1}~\mbox{on}~\mathbb{S}^{n-1}.$$
\end{enumerate}
\end{defn}




\begin{thm}\label{L-B}
Let $M$ be a positive constant and $n\geq3$.

\begin{enumerate}
\item[{\rm (a)}] For $m\geq2$, let
$f\in\mathcal{F}_{\varphi_{1},\ldots,\varphi_{m}}^{M}$.
Then there is a positive constant $R_{1}$ depending only on $M$ and $\varphi_{k}$ for $k\in\{1,\ldots,m\}$ such that
$\mathbb{B}^{n}(0,R_{1})\subset f(\mathbb{B}^{n})$.

\item[{\rm (b)}] For $m=1$, let $f\in\mathcal{F}_{\varphi_{1}}^{M}$.
Then $f(\mathbb{B}^{n})$ contains a ball
$\mathbb{B}^{n}(0,R_{1})$ with the radius $R_{1}$ satisfying
$$R_{1}\geq\frac{3n(1+3^{n-1})}{4}r_{1}M,$$
where
$r_{1}$ is an unique solution of the   equation: \beqq
&&\frac{1}{\left(\frac{nM}{2}+\frac{n\|\varphi_{1}\|_{\infty}}{n+1}\right)^{n-1}}-\frac{3n(1+3^{n-1})}{2}r_{1}M\\
&-&\|\varphi_{1}\|_{\infty} \max_{0\leq|x|\leq
r_{1}}\left(\int_{\mathbb{B}^{n}}\left|\nabla_{x}
G(x,y_{1})-\nabla_{x} G(0,y_{1})\right| dV(y_{1})\right)=0.\eeqq
\end{enumerate}
\end{thm}


A homeomorphism $f:~\Omega\rightarrow\Omega'$ between two open
subsets $\Omega$ and $\Omega'$ of $\mathbb{R}^{n}$ will be called a
$K$-quasiconformal mapping if

\begin{enumerate}
\item $f$ is an absolutely continuous function in almost every segment
parallel to some of the coordinate axes, and there exist the partial
derivatives which are locally $L^{n}$ integrable functions on
$\Omega$ (briefly, $f\in ACL^{n}$), and

\item $f$ satisfies the condition

\be\label{kqc}|D_{f}(x)|^{n}/K\leq J_{f}(x)\leq
K\big(l(D_{f}(x))\big)^{n}\ee
 at almost every $x$ in $\Omega$. 
\end{enumerate}
We remark that, for a continuous mapping $f$, the condition (1) is
equivalent to the condition that $f$ belongs to the Sobolev space
$W_{n,loc}^{1}(\Omega)$ (cf. \cite{Va,Vu}).

Given a subset $\Omega$ of $\mathbb{R}^{n}$, a function
$\psi:~\Omega\rightarrow\mathbb{R}^{n}$ is said to be {\it
bi-Lipschitz} if there is a  constant $c\geq1$ such that for all
 $x_{1},x_{2}\in\Omega$,

\be\label{eq-L-ck}\frac{1}{c}|x_{1}-x_{2}|\leq|\psi(x_{1})-\psi(x_{2})|\leq
c|x_{1}-x_{2}|.\ee Furthermore,  $\psi$ is called {\it Lipschitz} if
the right hand of (\ref{eq-L-ck}) holds, and $\psi$ is said to be
{\it co-Lipschitz} if it satisfies the left hand of (\ref{eq-L-ck}).

It is well known that  all sense-preserving bi-Lipschitz mappings
are quasiconformal mappings (cf. \cite{Va}). But quasiconformal
mappings are not necessarily bi-Lipschitz, not even Lipschitz (see
\cite{FV,KS,K-2011}). 


Pavlovi\'c \cite{Pav2} showed that harmonic quasiconformal mappings
of the unit disk $\mathbb{D}$ onto itself are bi-Lipschitz mappings.
In \cite{PS}, Partyka and Sakan improved Pavlovi\'c's corresponding
result and obtained an asymptotically sharp version. By using the
regularity theory of elliptic PDE's, Kalaj and Pavlovi\'c \cite{K1}
generalized the Lipschitz-property of harmonic quasiconformal
mappings to the quasiconformal solutions of Poisson's equations. The
same problem in the space is much more complicated because of the
lack of the techniques of complex analysis. It is well known that
the harmonic extension of a homeomorphism of the unit circle is
always a diffeomorphism of the unit disk $\mathbb{D}$. However, in
higher dimensions, the situation is quite different. Namely, Melas
\cite{Mel} constructed a homeomorphism of
$\mathbb{S}^{n-1}~(n\geq3)$ whose harmonic extension fails to be
diffeomorphic. 
On the discussion of the related topic, we refer to
\cite{BH-2004,CK-2019,CLW-2018,CW-2019,KS,K3,K-2011,KM-2010,
K-2008,MN, M-1968} and the  references therein. By using Theorem
\ref{thm-1} and Green's potential theory,
 we obtain the asymptotically sharp Lipschitz constant
which depends on the quasiconformal constant  $K$ and the Dirichlet
boundary value condition.


\begin{thm} \label{thm-2.1}
Let $K\geq1$, $n\geq3$, $m\geq2$,
$\varphi_{m}\in\mathcal{C}(\overline{\mathbb{B}^{n}},\mathbb{R}^{n})$
and $\varphi_{k}\in \mathcal{C}(\mathbb{S}^{n-1},\mathbb{R}^{n})$
for $k\in\{1,\ldots,m-1\}$.  Suppose that $f$ is a
$K$-quasiconformal self-mapping of $\mathbb{B}^{n}$ satisfying
$f(0)=0$ and \eqref{eq-ch-1} with the Dirichlet boundary value
condition:
$\Delta^{1}f=\varphi_{1},~\ldots,~\Delta^{m-1}f=\varphi_{m-1}$ on
$\mathbb{S}^{n-1}$. Then there are nonnegative constants
$N_{1}(K,\varphi_{1},\cdots,\varphi_{n})$ and $M_{1}(n,K)$
  with
$$\lim_{K\rightarrow1}M_{1}(n,K)=1~\mbox{and}~\lim_{\|\varphi_{1}\|_{\infty}\rightarrow0,\cdots,\|\varphi_{n}\|_{\infty}\rightarrow0}N_{1}(K,\varphi_{1},\cdots,\varphi_{n})=0$$
such that for all $x_{1}$ and $x_{2}$ in $\mathbb{B}^{n}$,

$$
|f(x_{1})-f(x_{2})|\leq
\big(M_{1}(n,K)+N_{1}(K,\varphi_{1},\cdots,\varphi_{n})\big)|x_{1}-x_{2}|.
$$
\end{thm}


We will give  several auxiliary results in Section
\ref{csw-sec2}. The proofs of Theorems \ref{thm-1} and \ref{thm-2}
will be presented in Section \ref{csw-sec2.0}, and the proof of
Theorem \ref{L-B} will be given in Section \ref{csw-sec3.0}. Theorem
\ref{thm-2.1} will be proved in the last Section.

\section{Auxiliary results}\label{csw-sec2}

\subsection{M\"obius Transformations of the Unit Ball}\label{sbcsw-sec2.1}

For $x\in\mathbb{B}^{n}$, the {\it M\"obius transformation} in
$\mathbb{B}^{n}$ is defined by \be\label{eq-ex1}
\phi_{x}(y)=\frac{|x-y|^{2}x-(1-|x|^{2})(y-x)}{[x,y]^{2}},~y\in\mathbb{B}^{n}.
\ee The set of isometries of the hyperbolic unit ball is a {\it
Kleinian subgroup} of all M\"obius transformations of the extended
spaces $\mathbb{R}^{n}\cup\{\infty\}$ onto itself. In the following,
we make use of the {\it automorphism group} $\Aut(\mathbb{B}^{n})$
consisting of all M\"obius transformations of the unit ball
$\mathbb{B}^{n}$ onto itself. We recall the following facts from
\cite{Bea}: For $x\in\mathbb{B}^{n}$ and
$\phi_{x}\in\Aut(\mathbb{B}^{n})$, we have $\phi_{x}(0)=x$,
$\phi_{x}(x)=0$, $\phi_{x}(\phi_{x}(y))=y \in\mathbb{B}^{n}$,
\be\label{II} |\phi_{x}(y)|=\frac{|x-y|}{[x,y]}\ee and
\be\label{III} |J_{\phi_{x}}(y)|=\frac{(1-|x|^{2})^{n}}{[x,y]^{2n}}.
\ee

\subsection{Gauss Hypergeometric Functions}\label{sbcsw-sec2.2}

For $a, b, c\in\mathbb{R}$ with $c\neq0, -1, -2, \ldots,$ the {\it
hypergeometric} function is defined by the power series in the variable $x$
$$F(a,b;c;x)=\sum_{k=0}^{\infty}\frac{(a)_{k}(b)_{k}}{(c)_{k}}\frac{x^{k}}{k!},~|x|<1.
$$
Here $(a)_{0}=1$,  $(a)_{k}=a(a+1)\cdots(a+k-1)$ for $k=1, 2, \ldots$, and
generally $(a)_{k}=\Gamma(a+k)/\Gamma(a)$ is the {\it Pochhammer} symbol,
where $\Gamma$ is the {\it Gamma function}. In particular, for $a, b, c>0$ and
$a+b<c$, we have (cf. \cite{PBM}) 
$$F(a,b;c;1)=\lim_{x\rightarrow1^{-}}
F(a,b;c;x)=\frac{\Gamma(c)\Gamma(c-a-b)}{\Gamma(c-a)\Gamma(c-b)}<\infty.
$$

The following result is useful in showing one of our main
results of the paper.

\begin{Lem}{\rm (\cite{K2}~$\mbox{or}$~\cite[2.5.16(43)]{PBM})}\label{pro-1}
For $\lambda_{1}>1$ and $\lambda_{2}>0$, we have
$$\int_{0}^{\pi}\frac{\sin^{\lambda_{1}-1}t}{(1+r^{2}-2r\cos t)^{\lambda_{2}}}dt=
\mathbf{B}\left(\frac{\lambda_{1}}{2},\frac{1}{2}\right)
F\big(\lambda_{2},\lambda_{2}+\frac{1-\lambda_{1}}{2};\frac{1+\lambda_{1}}{2};r^{2}\big),
$$
where $\mathbf{B}(.,.)$ denotes the beta function and $r\in[0,1)$.
\end{Lem}

\subsection{The spherical coordinates}\label{sbcsw-sec2.3}
Throughout this article, by $S$ and $T$ we denote the spherical
coordinates:

$$S:~Q_{0}^{n}=[0,1]\times[0,\pi]\times\cdots\times[0,\pi]\times[0,2\pi]\mapsto\mathbb{B}^{n}$$
and
$$T:~Q^{n-1}=[0,\pi]\times\cdots\times[0,\pi]\times[0,2\pi]\mapsto\mathbb{S}^{n-1},$$
$\big(S(r,\theta_{1},\cdots,\theta_{n-2},\theta_{n-1})=rT(\theta_{1},\cdots,\theta_{n-2},\theta_{n-1})\big)$,
defined by $S=(x_{1},x_{2},\cdots,x_{n-1})',$

\begin{eqnarray*}
x_{1}&=&r\cos\theta_{1},\\
x_{2}&=&r\sin\theta_{1}\sin\theta_{2},\\
&\vdots&\\
x_{n-1}&=&r\sin\theta_{1}\sin\theta_{2}\cdots\sin\theta_{n-2}\cos\theta_{n-1},\\
x_{n}&=&r\sin\theta_{1}\sin\theta_{2}\cdots\sin\theta_{n-2}\sin\theta_{n-1}.
\end{eqnarray*}
Then we have $$\det
D_{S}(r,\theta)=r^{n-1}\sin^{n-2}\theta_{1}\cdots\sin\theta_{n-2},$$
where $\theta=(\theta_{1},\ldots,\theta_{n-1})'$.



\section{The  heinz type inequalities for mappings satisfying polyharmonic equations}\label{csw-sec2.0}


The following result  easily follows from \cite[Theorem 1]{CP-Hi}.

\begin{Lem}\label{CP-1} Let $G(x,y)$ be the Green function defined in
{\rm (\ref{eq-ex0})}. Then for $x\in\mathbb{B}^{n}$,
$$\int_{\mathbb{B}^{n}}|G(x,y)|dV(y)=\frac{1-|x|^{2}}{2n}.$$
\end{Lem}

\begin{lem}\label{Ch-1.0}
Let $G(x,y)$ be the Green function defined in {\rm (\ref{eq-ex0})}.
Then for $x\in\mathbb{B}^{n}$,
$$\int_{\mathbb{B}^{n}}(1-|y|^{2})|G(x,y)|dV(y)=\frac{\big(n+4-n|x|^{2}\big)(1-|x|^{2})}{4n(n+2)}
\leq\frac{(n+4)(1-|x|^{2})}{4n(n+2)}.$$
\end{lem}
\bpf Let
$$I_{1}(x):=\int_{\mathbb{B}^{n}}(1-|y|^{2})|G(x,y)|dV(y).$$
For $x,y\in\mathbb{B}^{n}$ with $x\neq y$ and $|x|+|y|\neq0$, let
$z=\phi_{x}(y)$, where $\phi_{x}\in\Aut(\mathbb{B}^{n})$. Then
$y=\phi_{x}(z)$ and
\be\label{eq-4d}1-|\phi_{x}(z)|^{2}=\frac{(1-|x|^{2})(1-|z|^{2})}{[x,z]^{2}}.\ee
It follows from \eqref{eq-ex1}  that

\begin{eqnarray*}
x-\phi_{x}(z)=\frac{x[x,z]^{2}-|x-z|^{2}x+(1-|x|^{2})(z-x)}{[x,z]^{2}}
=\frac{(z-x|z|^{2})(1-|x|^{2})}{[x,z]^{2}},
\end{eqnarray*}
which gives \be\label{eq-3}
|x-\phi_{x}(z)|=\frac{|z|(1-|x|^{2})}{[x,z]}. \ee By \eqref{II}, we
have
\begin{eqnarray*}
\left|\frac{1}{|x-y|^{n-2}}-\frac{1}{[x,y]^{n-2}}\right|&=&\frac{1}{|x-y|^{n-2}}\left|1-\frac{|x-y|^{n-2}}{[x,y]^{n-2}}\right|\\
\nonumber&=&\frac{1-|z|^{n-2}}{|x-\phi_{x}(z)|^{n-2}},
\end{eqnarray*} which, together with (\ref{eq-3}), implies that
 \be\label{eq-4}
\left|\frac{1}{|x-y|^{n-2}}-\frac{1}{[x,y]^{n-2}}\right|=\frac{[x,z]^{n-2}(1-|z|^{n-2})}{|z|^{n-2}(1-|x|^{2})^{n-2}}.
\ee

Using the spherical coordinates and Lemma A, we obtain
\beq\label{eq-6f}
\int_{\mathbb{S}^{n-1}}\frac{d\sigma(\zeta)}{|rx-\zeta|^{4+n}}&=&\frac{1}{\int_{0}^{\pi}\sin^{n-2}t~dt}\int_{0}^{\pi}
\frac{\sin^{n-2}t}{\left(1+r^{2}|x|^{2}-2r|x|\cos
t\right)^{\frac{n+4}{2}}}dt\\ \nonumber
&=&\frac{\Gamma\big(\frac{n}{2}\big)}{\sqrt{\pi}\Gamma\big(\frac{n-1}{2}\big)}\cdot
\frac{\sqrt{\pi}\Gamma\big(\frac{n-1}{2}\big)}{\Gamma\big(\frac{n}{2}\big)}F\Big(\frac{n+4}{2},3;\frac{n}{2};r^{2}|x|^{2}\Big)\\
 \nonumber&=& F\Big(\frac{n+4}{2},3;\frac{n}{2};r^{2}|x|^{2}\Big)\\
\nonumber
&=&\sum_{k=0}^{\infty}\frac{(k+1)(k+2)(n+2k)(n+2k+2)}{2n(n+2)}r^{2k}|x|^{2k}.
\eeq By (\ref{eq-4d}), (\ref{eq-3}), (\ref{eq-4}), (\ref{eq-6f}) and
the change of variables, we obtain

\beqq\label{eq-ch-2g} I_{1}(x)&=&\int_{\mathbb{B}^{n}}
\frac{(1-|z|^{n-2})(1-|x|^{2})^{3}(1-|z|^{2})}{|z|^{n-2}[x,z]^{n+4}}dV(z)\\
\nonumber
&=&\frac{(1-|x|^{2})^{3}}{n-2}\int_{0}^{1}\left[r(1-r^{2})(1-r^{n-2})\int_{\partial\mathbb{B}^{n}}\frac{d\sigma(\zeta)}{|rx-\zeta|^{n+4}}\right]dr\\
\nonumber
&=&\frac{(1-|x|^{2})^{3}}{n-2}\sum_{k=0}^{\infty}\frac{(k+1)(k+2)(n+2k)(n+2k+2)}{2n(n+2)}|x|^{2k}\\
\nonumber &&\times\int_{0}^{1}r^{2k+1}(1-r^{2})(1-r^{n-2})dr\\
\nonumber
&=&\frac{\big(n+4-n|x|^{2}\big)(1-|x|^{2})}{4n(n+2)}. \eeqq The
proof of this lemma is complete. \epf

\subsection*{Proof of Theorem \ref{thm-1}} By (\ref{eq-ch-3.0}), we have

$$f(x)=P[\varphi_{0}](x)+\sum_{k=1}^{m}(-1)^{k}G_{k}[\varphi_{k}](x),~x\in\mathbb{B}^{n},$$
where $\varphi_{0}=f|_{\mathbb{S}^{n-1}}$, $G_{k}[\varphi_{k}]$ are
defined in (\ref{eq-ch-3.1}) for $k\in\{ 1,\ldots,m-1\}$, and
$G_{m}[\varphi_{m}]$ is defined in (\ref{eq-ch-3.2}). Next, we
estimate $|G_{m}[\varphi_{m}]|$ and
$|G_{k}[\varphi_{k}]|$ for $k\in\{1,\ldots,m-1\}$. 


\noindent $\mathbf{Case~ 1.}$  $m=2$ and $k=1$.

By Lemma B, we have

\beq\label{eq-ch-3g} |G_{1}[\varphi_{1}](x)|&=&
\left|\int_{\mathbb{B}^{n}}G(x,y_{1})\left(\int_{\mathbb{S}^{n-1}}P(y_{1},\zeta)\varphi_{1}(\zeta)d\sigma(\zeta)\right)dV(y_{1})\right|\\
\nonumber &\leq&\frac{\|\varphi_{1}\|_{\infty}}{2n}(1-|x|^{2}).\eeq

\noindent $\mathbf{Case~ 2.}$ $m\geq3$ and $2\leq k\leq m-1$.


It follows from (\ref{eq-ch-3g}) and Lemma \ref{Ch-1.0} that

\beq\label{eq-ch-4g} |G_{k}[\varphi_{k}](x)|&=&
\bigg|\int_{\mathbb{B}^{n}}\cdots\int_{\mathbb{B}^{n}}G(x,y_{1})\cdots
G(y_{k-1},y_{k})\\ \nonumber
&&\times\left(\int_{\mathbb{S}^{n-1}}P(y_{k},\zeta)\varphi_{k}(\zeta)d\sigma(\zeta)\right)dV(y_{k})\cdots
dV(y_{1})\bigg|\\ \nonumber
&\leq&\frac{\|\varphi_{k}\|_{\infty}}{2n}\int_{\mathbb{B}^{n}}\cdots\int_{\mathbb{B}^{n}}|G(x,y_{1})|\cdots
|G(y_{k-2},y_{k-1})|\\
\nonumber&&\times(1-|y_{k-1}|^{2})dV(y_{k-1})\cdots dV(y_{1})\,\,\,\mbox{(by (\ref{eq-ch-3g}))}\\
\nonumber
&\leq&\frac{\|\varphi_{k}\|_{\infty}}{2n}\left[\frac{n+4}{4n(n+2)}\right]^{k-1}(1-|x|^{2}).\,\,\,\mbox{(by
Lemma \ref{Ch-1.0})} \eeq

Now we estimate $|G_{m}[\varphi_{m}]|$, where  $m\geq2$.

By (\ref{eq-ch-3g}) and Lemma \ref{Ch-1.0}, we see that

\beq\label{eq-ch-6g}
|G_{m}[\varphi_{m}](x)|&\leq&\frac{\|\varphi_{m}\|_{\infty}}{2n}
\int_{\mathbb{B}^{n}}\cdots\int_{\mathbb{B}^{n}}|G(x,y_{1})|\cdots |G(y_{m-2},y_{m-1})|\\
\nonumber &&\times(1-|y_{m-1}|^{2})dV(y_{m-1})\cdots dV(y_{1})\,\,\,\mbox{(by (\ref{eq-ch-3g}))}\\
\nonumber
&\leq&\frac{\|\varphi_{m}\|_{\infty}}{2n}\left[\frac{n+4}{4n(n+2)}\right]^{m-1}(1-|x|^{2}).
\,\,\,\mbox{(by Lemma \ref{Ch-1.0})} \eeq

Therefore, it follows from (\ref{eq-K}), (\ref{eq-ch-3g}),
(\ref{eq-ch-4g}) and (\ref{eq-ch-6g}) that

\begin{eqnarray*}
\left|f(x)-\frac{1-|x|^{2}}{(1+|x|^{2})^{\frac{n}{2}}}P[\varphi_{0}](0)\right|
&\leq&\left|P[\varphi_{0}](x)-\frac{1-|x|^{2}}{(1+|x|^{2})^{\frac{n}{2}}}P[\varphi_{0}](0)\right|+\sum_{k=1}^{m}|G_{k}[\varphi_{k}](x)|\\
&\leq&\|\varphi_{0}\|_{\infty}U(|x|N)\\
&&+\sum_{k=1}^{m}\frac{\|\varphi_{k}\|_{\infty}}{2n}\left[\frac{n+4}{4n(n+2)}\right]^{k-1}(1-|x|^{2}).
\end{eqnarray*} The proof of the theorem is complete. \qed

\begin{Lem}\label{K-5}{\rm (\cite[Lemma 2.3]{K5})} For $r\in[0,1]$, let $\Phi(r)=\partial U(rN)/\partial r.$ Then
$\Phi(r)$ is decreasing on $r\in[0,1]$ and

$$\label{tt-1}\Phi(r)\geq\frac{\partial U(rN)}{\partial r}\bigg|_{r=1}=
\frac{n!\big[1+n-(n-2)F\big(\frac{1}{2},1;\frac{n+3}{2};-1\big)\big]}{2^{\frac{3n}{2}}\Gamma(\frac{n+1}{2})\Gamma(\frac{n+3}{2})}.$$

\end{Lem}

\subsection*{Proof of Theorem \ref{thm-2}}
By (\ref{eq-ch-3.0}), we have

$$f(x)=P[\varphi_{0}](x)+\sum_{k=1}^{m}(-1)^{k}G_{k}[\varphi_{k}](x),~x\in\mathbb{B}^{n},$$
where   $G_{k}[\varphi_{k}]$ are defined in (\ref{eq-ch-3.1})  for
 $k\in\{1,\ldots,m-1\}$, and $|G_{m}[\varphi_{m}]|$  is
defined in (\ref{eq-ch-3.2}). By the assumption, we see that
$$0=f(0)=P[\varphi_{0}](0)+\sum_{k=1}^{m}(-1)^{k}G_{k}[\varphi_{k}](0),$$
which, together with (\ref{eq-ch-4g}),  (\ref{eq-ch-6g}) and Theorem
\ref{thm-1}, implies that

\beq\label{he-1}
|f(\zeta)-f(r\zeta)|&=&\bigg|f(\zeta)+\frac{1-r^{2}}{(1+r^{2})^{\frac{n}{2}}}P[\varphi_{0}](0)\\
\nonumber
&&+\sum_{k=1}^{m}(-1)^{k}G_{k}[\varphi_{k}](0)\frac{1-r^{2}}{(1+r^{2})^{\frac{n}{2}}}-f(r\zeta)\bigg|\\
\nonumber
&\geq&1-\left|f(r\zeta)-\frac{1-r^{2}}{(1+r^{2})^{\frac{n}{2}}}P[\varphi_{0}](0)\right|\\
\nonumber&&-\frac{1-r^{2}}{(1+r^{2})^{\frac{n}{2}}}
\sum_{k=1}^{m}|G_{k}[\varphi_{k}](0)|\\ \nonumber&\geq&
1-U(rN)-\sum_{k=1}^{m}\frac{\|\varphi_{k}\|_{\infty}}{2n}\left[\frac{n+4}{4n(n+2)}\right]^{k-1}(1-r^{2})\\
\nonumber
&&-\sum_{k=1}^{m}\frac{\|\varphi_{k}\|_{\infty}}{2n}\left[\frac{n+4}{4n(n+2)}\right]^{k-1}\frac{1-r^{2}}{(1+r^{2})^{\frac{n}{2}}},
 \eeq where $r\in[0,1)$.
On the other hand, for $x\in\mathbb{B}^{n}$, there is a
$\rho\in(|x|,1)$ such that

\be\label{he-2} \frac{1-U(rN)}{1-r}=\frac{\partial U(\rho
N)}{\partial r},\ee where $r=|x|$. It follows from (\ref{he-1}),
(\ref{he-2}) and Lemma C that

\begin{eqnarray*}
\liminf_{r\rightarrow1^{-}}\frac{|f(\zeta)-f(r\zeta)|}{1-r}&\geq&\liminf_{r\rightarrow1^{-}}\frac{1-U(rN)}{1-r}
-\sum_{k=1}^{m}\frac{\|\varphi_{k}\|_{\infty}}{n}\left[\frac{n+4}{4n(n+2)}\right]^{k-1}\\
\nonumber
&&-\sum_{k=1}^{m}\frac{\|\varphi_{k}\|_{\infty}}{n}\left[\frac{n+4}{4n(n+2)}\right]^{k-1}\frac{1}{2^{\frac{n}{2}}}\\
&\geq&\frac{n!\big[1+n-(n-2)F\big(\frac{1}{2},1;\frac{n+3}{2};-1\big)\big]}{2^{\frac{3n}{2}}\Gamma(\frac{n+1}{2})\Gamma(\frac{n+3}{2})}\\
&&-\sum_{k=1}^{m}\frac{\|\varphi_{k}\|_{\infty}}{n}\left[\frac{n+4}{4n(n+2)}\right]^{k-1}\left(1+\frac{1}{2^{\frac{n}{2}}}\right).
\end{eqnarray*}

At last, we prove the sharpness part. Especially, if
$\|\varphi_{k}\|_{\infty}=0$ for $k\in\{1,\ldots,m\}$, then the
sharpness part easily follows from \cite[Theorem 2.5]{K5}. The proof
of this theorem is complete. \qed

\section{ Bloch type theorem for mappings  satisfying  polyharmonic equations}\label{csw-sec3.0}




The main purpose of this section is to prove Theorem \ref{L-B}. We
start with some lemmas which are useful  to the proof of Theorem \ref{L-B}.

\begin{Thm}{\rm (see \cite[Theorem 2.7]{KMP})}\label{Result-2022}
For $n\geq3$, let $u$ be a harmonic function of $\mathbb{B}^{n}$ into itself. Then
$$|D_{u}(x)|\leq\frac{n}{2}\frac{1}{1-|x|}.$$
\end{Thm}

A matrix-valued function $A(x)=\big(a_{ij}(x)\big)_{n\times n}$is
called {\it matrix-valued  harmonic function} if each of its entries
$a_{ij}(x)$ is a  harmonic function from an open subset
$\Omega\subset\mathbb{R}^{n}$ into $\mathbb{R}$.

\begin{Lem}\label{CMPW-Lem}{\rm (\cite[Lemma 3.1]{CMPW})} For $r>0$,
let $A(x)=\big(a_{ij}(x)\big)_{n\times n}$ be a matrix-valued
harmonic function defined on  $\mathbb{B}^{n}(0,r)$. If $A(0)=0$ and
$|A(x)|\leq M$ in $\mathbb{B}^{n}(0,r),$ then
$$|A(x)|\leq M\left(1-\frac{r^{n-2}(r-|x|)}{(r+|x|)^{n-1}}\right),
$$ where $M$ is a positive constant.
\end{Lem}

\begin{Lem}\label{K-3}{\rm (\cite[Lemma 2.5]{K-2011}~\mbox{or}~\cite[p. 24-26]{Vl})}
Let $\varrho$ be a bounded (absolutely) integrable function defined
on a bounded domain $\Omega\subset\mathbb{R}^{n}$. Then the
potential type integral

$$\Psi(x)=\int_{\Omega}\frac{\varrho(y)dV(y)}{|x-y|^{\alpha}}$$ belongs to
the space $\mathcal{C}^{k}(\mathbb{R}^{n},\mathbb{R})$, where
$k+\alpha<n.$ Moreover,

$$\nabla\Psi(x)=\int_{\Omega}\nabla\left(\frac{1}{|x-y|^{\alpha}}\varrho(y)\right)dV(y).$$
\end{Lem}

\begin{lem}\label{eq-jj-1} Suppose that $\varphi_{k}\in
\mathcal{C}(\mathbb{S}^{n-1},\mathbb{R}^{n})$  and
$G_{k}[\varphi_{k}]$ are defined in {\rm (\ref{eq-ch-3.1})}, where
 $m\in\mathbb{N}\setminus\{1\}$, $k\in\{1,\ldots,m-1\}$ and
 $n\geq3.$ Then, for
$x\in\mathbb{B}^{n}$,

\be\label{eq-ff-1}\left|D_{G_{k}[\varphi_{k}]}(x)\right|\leq\begin{cases}
\displaystyle \frac{n}{n+1}\|\varphi_{1}\|_{\infty},
&\mbox{ if }\, k=1,\\
\displaystyle
\frac{\|\varphi_{k}\|_{\infty}}{2n}\left[\frac{n+4}{4n(n+2)}\right]^{k-2}\delta(n),
&\mbox{ if }\, 2\leq k\leq m-1,
\end{cases}\ee where $$\delta(n)=\frac{(n^{2}-4)}{3(n^{2}-1)}c_{0}+\frac{4}{n(n+2)}~\mbox{and}~c_{0}=\max_{0\leq t<1}[(2-t^{2})(1+t)]\approx2.631.$$
Moreover, $D_{G_{k}[\varphi_{k}]}$ has a continuous extension to the
boundary and, for $\varepsilon\in\mathbb{S}^{n-1}$,

\be\label{eq-ff-2}\left|D_{G_{k}[\varphi_{k}]}(\varepsilon)\right|\leq\begin{cases}
\displaystyle \frac{1}{n}\|\varphi_{1}\|_{\infty},
&\mbox{ if }\, k=1,\\
\displaystyle
\frac{\|\varphi_{k}\|_{\infty}}{2n^{2}(n+2)}\left[\frac{n+4}{4n(n+2)}\right]^{k-2},
&\mbox{ if }\, 2\leq k\leq m-1.
\end{cases}\ee
\end{lem}

\bpf We divide the proof of this lemma into two steps.
\bst\label{t1} We first prove {\rm (\ref{eq-ff-2})}. \est 

 \noindent $\mathbf{Case~ 1.}$
$m\geq3$ and $2\leq k\leq m-1$.

For $k\in\{2,\ldots,m-1\}$ and $m\geq3$, let \beq\label{eq-hh-1}
 \varrho_{\varphi_{k}}(y_{1})&=&
\int_{\mathbb{B}^{n}}\bigg(G(y_{1},y_{2})\cdots\int_{\mathbb{B}^{n}}\bigg(G(y_{k-1},y_{k})\\
\nonumber
&&\times\left(\int_{\mathbb{S}^{n-1}}P(y_{k},\zeta)\varphi_{k}(\zeta)d\sigma(\zeta)\right)dV(y_{k})\bigg)\cdots\bigg)
dV(y_{2}).\eeq Then, by (\ref{eq-ch-4g}), we have

\be\label{eq-hh2}|\varrho_{\varphi_{k}}(y_{1})|\leq\frac{\|\varphi_{k}\|_{\infty}}{2n}\left[\frac{n+4}{4n(n+2)}\right]^{k-2}(1-|y_{1}|^{2})
<\infty,\ee where $y_{1}\in\mathbb{B}^{n}.$

Applying Lemma F to $$\int_{\mathbb{B}^{n}}G(x,y_{1})
\varrho_{\varphi_{k}}(y_{1})dV(y_{1}),$$ we see that, for any fixed
$\xi\in\mathbb{R}^{n}$,

\be\label{eq-hh4}
D_{G_{k}[\varphi_{k}]}(x)\xi=\int_{\mathbb{B}^{n}}\langle\nabla_{x}
G(x,y_{1}),\xi\rangle\varrho_{\varphi_{k}}(y_{1})dV(y_{1}).\ee

By calculations, we obtain

$$\nabla_{x} G(x,y_{1})=-\frac{1}{\omega_{n-1}}\left(\frac{x-y_{1}}{|x-y_{1}|^{n}}-\frac{|y_{1}|^{2}x-y_{1}}{[x,y_{1}]^{n}}\right),$$
which gives \be\label{eq-hh3}\left|\nabla_{x}
G(x,y_{1})\right|\leq\frac{1}{\omega_{n-1}}\left(\frac{1}{|x-y_{1}|^{n-1}}+\frac{1}{\big||y_{1}|^{2}x-y_{1}\big|^{n-1}}\right),\ee
where $x\in\mathbb{B}^{n}.$ It follows from (\ref{eq-hh2}),
(\ref{eq-hh4}), (\ref{eq-hh3}), Lemma F and the Lebesgue
Dominated Convergence Theorem  that, for
$\varepsilon\in\mathbb{S}^{n-1}$,

\beq\label{eq-hh5}
D_{G_{k}[\varphi_{k}]}(\varepsilon)\xi&=&\lim_{x\rightarrow\varepsilon}\int_{\mathbb{B}^{n}}\langle\nabla_{x}
G(x,y_{1}),\xi\rangle\varrho_{\varphi_{k}}(y_{1})dV(y_{1})\\
\nonumber
&=&\int_{\mathbb{B}^{n}}\lim_{x\rightarrow\varepsilon}\langle\nabla_{x}
G(x,y_{1}),\xi\rangle\varrho_{\varphi_{k}}(y_{1})dV(y_{1})\\
\nonumber
&=&\frac{1}{\omega_{n-1}}\int_{\mathbb{B}^{n}}\langle\varepsilon,\xi\rangle
\frac{1-|y_{1}|^{2}}{|\varepsilon-y_{1}|^{n}}\varrho_{\varphi_{k}}(y_{1})dV(y_{1}).
\eeq

Next, we estimate $|D_{G_{k}[\varphi_{k}]}(\varepsilon)|.$ By
(\ref{eq-hh2}) and (\ref{eq-hh5}), we have

\begin{eqnarray*}
|D_{G_{k}[\varphi_{k}]}(\varepsilon)\xi|&\leq&\frac{\|\varphi_{k}\|_{\infty}}{2n}\left[\frac{n+4}{4n(n+2)}\right]^{k-2}\frac{|\xi|}{\omega_{n-1}}\int_{\mathbb{B}^{n}}
\frac{(1-|y_{1}|^{2})^{2}}{|\varepsilon-y_{1}|^{n}}dV(y_{1})\\
&=&\frac{\|\varphi_{k}\|_{\infty}}{2n^{2}(n+2)}\left[\frac{n+4}{4n(n+2)}\right]^{k-2}|\xi|,
\end{eqnarray*} which implies that
$$|D_{G_{k}[\varphi_{k}]}(\varepsilon)|\leq\frac{\|\varphi_{k}\|_{\infty}}{2n^{2}(n+2)}\left[\frac{n+4}{4n(n+2)}\right]^{k-2}.$$

\noindent $\mathbf{Case~ 2.}$  $m=2$ and $k=1$.

Using the spherical coordinates and Proposition \ref{pro-1}, we
obtain
\beqq
\left|D_{G_{1}[\varphi_{1}]}(\varepsilon)\xi\right|
&=&\bigg|\frac{1}{\omega_{n-1}}\int_{\mathbb{B}^{n}}\langle\varepsilon,\xi\rangle
\frac{1-|y_{1}|^{2}}{|\varepsilon-y_{1}|^{n}}\left(\int_{\mathbb{S}^{n-1}}P(y_{1},\zeta)\varphi_{1}(\zeta)d\sigma(\zeta)\right)dV(y_{1})\bigg|\\
&\leq&\frac{\|\varphi_{1}\|_{\infty}|\xi|}{\omega_{n-1}}\int_{\mathbb{B}^{n}}\frac{1-|y_{1}|^{2}}{|\varepsilon-y_{1}|^{n}}dV(y_{1})=
\frac{\|\varphi_{1}\|_{\infty}|\xi|}{n}, \eeqq which yields that

$$\left|D_{G_{1}[\varphi_{1}]}(\varepsilon)\right|\leq\frac{\|\varphi_{1}\|_{\infty}}{n}.$$

\bst\label{t2} Next, we show that {\rm (\ref{eq-ff-1})}. \est
\noindent $\mathbf{Case~ 3.}$ $m\geq3$ and $2\leq k\leq m-1$.

In order to estimate $|D_{G_{k}[\varphi_{k}]}(x)\xi|,$ we first show
that, for $x\in\mathbb{B}^{n},$ \be\label{hj-1}I_{2}(x)\leq
\delta(n),\ee where
$$I_{2}(x)=\int_{\mathbb{B}^{n}}\left|\nabla_{x}
G(x,y_{1})\right|(1-|y_{1}|^{2})dV(y_{1}).$$

In order to prove (\ref{hj-1}), we let $z=\phi_{x}(y_{1})$, where
$\phi_{x}\in\Aut(\mathbb{B}^{n})$. Then, by (\ref{eq-4d}) and
(\ref{eq-3}), we have

\beqq \left|\nabla_{x}
G(x,y_{1})\right|(1-|y_{1}|^{2})&=&\frac{(1-|y_{1}|^{2})}{\omega_{n-1}}\left|\frac{x-y_{1}}{|x-y_{1}|^{n}}-\frac{|y_{1}|^{2}x-y_{1}}{[x,y_{1}]^{n}}\right|\\
&=&\frac{1}{\omega_{n-1}}\frac{(1-|y_{1}|^{2})}{|x-y_{1}|^{n}}\left|(x-y_{1})-\frac{\big(|y_{1}|^{2}x-y_{1}\big)|x-y_{1}|^{n}}{[x,y_{1}]^{n}}\right|\\
&=&\frac{\big|(x-\phi_{x}(z))-(|\phi_{x}(z)|^{2}x-\phi_{x}(z))|z|^{n}\big|}{\omega_{n-1}}\frac{(1-|\phi_{x}(z)|^{2})}{|x-\phi_{x}(z)|^{n}}\\
&=&\frac{1}{\omega_{n-1}}\frac{[x,z]^{n-4}(1-|z|^{2})}{|z|^{n-1}(1-|x|^{2})^{n-2}}\left|\frac{z}{|z|}-x|z|-z|z|^{n-1}+x|z|^{n-1}\right|,
\eeqq which, together with (\ref{III}), gives that

\beqq
I_{2}(x)&=&(1-|x|^{2})^{2}\int_{0}^{1}\left[(1-r^{2})\int_{\mathbb{S}^{n-1}}\frac{\big|\zeta-rx-r^{n}\zeta+xr^{n-1}\big|}{|xr-\zeta|^{n+4}}d\sigma(\zeta)\right]dr\\
&\leq&(1-|x|^{2})^{2}\int_{0}^{1}\left[(1-r^{2})\int_{\mathbb{S}^{n-1}}\frac{|rx-\zeta|(1-r^{n-2})+r^{n-2}(1-r^{2})}{|xr-\zeta|^{n+4}}d\sigma(\zeta)\right]dr\\
&=&I_{3}(x)+I_{4}(x), \eeqq where

$$I_{3}(x)=(1-|x|^{2})^{2}\int_{0}^{1}\left[(1-r^{2})(1-r^{n-2})\int_{\mathbb{S}^{n-1}}\frac{1}{|xr-\zeta|^{n+3}}d\sigma(\zeta)\right]dr$$
and
$$I_{4}(x)=(1-|x|^{2})^{2}\int_{0}^{1}\left[(1-r^{2})^{2}r^{n-2}\int_{\mathbb{S}^{n-1}}\frac{1}{|xr-\zeta|^{n+4}}d\sigma(\zeta)\right]dr.$$
By computations, we obtain
$$\int_{0}^{1}(1-r^{2})(1-r^{n-2})r^{2k}dr=\frac{2(n-2)(n+4k+2)}{(2k+1)(2k+3)(n+2k-1)(n+2k+1)}$$
and
$$\frac{(k+1)(n+2k)(n+2k+2)(n+4k+2)}{(2k+1)(2k+3)(n+2k-1)(n+2k+1)}(k+2)\leq\frac{n(n+2)^{2}}{3(n^{2}-1)}(k+2),$$
which, together with (\ref{eq-6f}), imply that

\beq\label{eq-s-1.1}\nonumber
I_{3}(x)&\leq&(1+|x|)(1-|x|^{2})^{2}\int_{0}^{1}\left[(1-r^{2})(1-r^{n-2})\int_{\mathbb{S}^{n-1}}\frac{1}{|xr-\zeta|^{n+4}}d\sigma(\zeta)\right]dr\\
\nonumber
&=&(1+|x|)(1-|x|^{2})^{2}\sum_{k=0}^{\infty}\frac{(k+1)(k+2)(n+2k)(n+2k+2)}{2n(n+2)}|x|^{2k}\\
\nonumber &&\times\int_{0}^{1}(1-r^{2})(1-r^{n-2})r^{2k}dr\\
&\leq&(1+|x|)(1-|x|^{2})^{2}\frac{(n^{2}-4)}{3(n^{2}-1)}\sum_{k=0}^{\infty}(k+2)|x|^{2k}
=\frac{(n^{2}-4)}{3(n^{2}-1)}c_{0},\eeq where
$$c_{0}=\max_{0\leq|x|<1}(2-|x|^{2})(1+|x|)\approx2.631.$$

It follows from (\ref{eq-6f}) and the inequality
$$\frac{(k+1)(k+2)(n+2k)(n+2k+2)}{(n+2k-1)(n+2k+3)(n+2k+1)}\leq\frac{k+2}{2}$$
that \beq\label{eq-s-1.2}\nonumber
\frac{I_{4}(x)}{(1-|x|^{2})^{2}}&=&\sum_{k=0}^{\infty}\frac{(k+1)(k+2)(n+2k)(n+2k+2)}{2n(n+2)}|x|^{2k}\int_{0}^{1}(1-r^{2})^{2}r^{n-2+2k}dr\\
\nonumber
&=&\sum_{k=0}^{\infty}\frac{4(k+1)(k+2)(n+2k)(n+2k+2)}{n(n+2)(n+2k-1)(n+2k+3)(n+2k+1)}|x|^{2k}\\
 &\leq&\frac{2}{n(n+2)}\sum_{k=0}^{\infty}
(k+2)|x|^{2k}=\frac{2(2-|x|^{2})}{n(n+2)(1-|x|^{2})^{2}}.\eeq Then
combining (\ref{eq-s-1.1}) and (\ref{eq-s-1.2}) gives the estimate
for (\ref{hj-1}). Therefore, we conclude from (\ref{eq-hh2}),
(\ref{eq-hh4}) and (\ref{hj-1}) that


\beqq
|D_{G_{k}[\varphi_{k}]}(x)\xi|&\leq&|\xi|\int_{\mathbb{B}^{n}}|\nabla_{x}
G(x,y_{1})||\varrho_{\varphi_{k}}(y_{1})|dV(y_{1})\\ \nonumber
&\leq&\frac{\|\varphi_{k}\|_{\infty}}{2n}\left[\frac{n+4}{4n(n+2)}\right]^{k-2}I_{2}(x)|\xi|\\
\nonumber
&\leq&\frac{\|\varphi_{k}\|_{\infty}}{2n}\left[\frac{n+4}{4n(n+2)}\right]^{k-2}\delta(n)|\xi|,
 \eeqq which yields that

$$|D_{G_{k}[\varphi_{k}]}(x)|\leq\frac{\|\varphi_{k}\|_{\infty}}{2n}\left[\frac{n+4}{4n(n+2)}\right]^{k-2}\delta(n).$$

\noindent $\mathbf{Case~ 4.}$  $m=2$ and $k=1$.

By \cite[Theorem 2.1]{K2},
\beqq
|D_{G_{1}[\varphi_{1}]}(x)\xi|&\leq&|\xi|\|\varphi_{1}\|_{\infty}\int_{\mathbb{B}^{n}}|\nabla_{x}
G(x,y_{1})|dV(y_{1})\leq\frac{n\|\varphi_{1}\|_{\infty}}{n+1}|\xi|,
\eeqq which implies that
$$|D_{G_{1}[\varphi_{1}]}(x)|\leq\frac{n\|\varphi_{1}\|_{\infty}}{n+1}.$$
The proof of this lemma is finished. \epf

\begin{lem}\label{eq-jj-2} Suppose that $\varphi_{m}\in
\mathcal{C}(\mathbb{B}^{n},\mathbb{R}^{n})$  and
$G_{m}[\varphi_{m}]$ are defined in {\rm (\ref{eq-ch-3.2})}, where
 $m\in\mathbb{N}\setminus\{1\}$ and
 $n\geq3.$ Then, for
$x\in\mathbb{B}^{n}$,

$$\left|D_{G_{m}[\varphi_{m}]}(x)\right|\leq\frac{\|\varphi_{m}\|_{\infty}}{2n}\left[\frac{n+4}{4n(n+2)}\right]^{m-2}\delta(n),$$
where $\delta(n)$ is defined in Lemma {\rm \ref{eq-jj-1}}.

Moreover, $D_{G_{m}[\varphi_{m}]}$ has a continuous extension to the
boundary and, for $\varepsilon\in\mathbb{S}^{n-1}$,
$$\left|D_{G_{m}[\varphi_{m}]}(\varepsilon)\right|\leq\frac{\|\varphi_{m}\|_{\infty}}{2n^{2}(n+2)}\left[\frac{n+4}{4n(n+2)}\right]^{m-2}.$$
\end{lem}

\bpf Let \beq \label{eq-jj-2.1}\nonumber
\varrho_{\varphi_{m}}(y_{1})&=&
\int_{\mathbb{B}^{n}}\bigg(G(y_{1},y_{2})\cdots\int_{\mathbb{B}^{n}}\bigg(G(y_{m-2},y_{m-1})\\
&&\times\left(\int_{\mathbb{B}^{n}}G(y_{m-1},y_{m})\varphi_{m}(y_{m})dV(y_{m})\right)dV(y_{m-1})\bigg)\cdots\bigg)
dV(y_{2}).\eeq By Lemma B and (\ref{eq-ch-2g}), we have

 \beq \label{eq-jj-3}
|\varrho_{\varphi_{m}}(y_{1})|&\leq&\|\varphi_{m}\|_{\infty}\int_{\mathbb{B}^{n}}\bigg(|G(y_{1},y_{2})|\cdots\\
\nonumber
&&\times\int_{\mathbb{B}^{n}}|G(y_{m-2},y_{m-1})|dV(y_{m-1})\cdots\bigg)
dV(y_{2})\\
\nonumber&\leq&\frac{\|\varphi_{m}\|_{\infty}}{2n}\left[\frac{(n+4)}{4n(n+2)}\right]^{m-2}(1-|y_{1}|^{2})\\
\nonumber&\leq&\frac{\|\varphi_{m}\|_{\infty}}{2n}\left[\frac{(n+4)}{4n(n+2)}\right]^{m-2},
 \eeq
which, together with (\ref{eq-hh3}), Lemma F and the
Lebesgue Dominated Convergence Theorem, implies that, for any fixed
$\xi\in\mathbb{R}^{n}$,

\be\label{eq-jj4}
D_{G_{m}[\varphi_{m}]}(x)\xi=\int_{\mathbb{B}^{n}}\langle\nabla_{x}
G(x,y_{1}),\xi\rangle\varrho_{\varphi_{m}}(y_{1})dV(y_{1}) \ee and
$D_{G_{m}[\varphi_{m}]}$ has a continuous extension to the boundary.

Next we estimate $|D_{G_{m}[\varphi_{m}]}(x)|$ for
$x\in\mathbb{B}^{n}$, and $|D_{G_{m}[\varphi_{m}]}(\eta)|$ for
$\eta\in\mathbb{S}^{n-1}$, respectively.

It follows from (\ref{hj-1}), (\ref{eq-jj-3}) and (\ref{eq-jj4})
that

\begin{eqnarray*}
|D_{G_{m}[\varphi_{m}]}(x)\xi|&\leq&\int_{\mathbb{B}^{n}}\big|\langle\nabla_{x}
G(x,y_{1}),\xi\rangle\big||\varrho_{\varphi_{m}}(y_{1})|dV(y_{1})\\&\leq&\frac{\|\varphi_{m}\|_{\infty}}{2n}\left[\frac{n+4}{4n(n+2)}\right]^{m-2}\delta(n)|\xi|,
\end{eqnarray*} and
\begin{eqnarray*}
|D_{G_{m}[\varphi_{m}]}(\varepsilon)\xi|&=&\left|\lim_{x\rightarrow\varepsilon}\int_{\mathbb{B}^{n}}\langle\nabla_{x}
G(x,y_{1}),\xi\rangle\varrho_{\varphi_{m}}(y_{1})dV(y_{1})\right|\\
&=&\frac{1}{\omega_{n-1}}\int_{\mathbb{B}^{n}}\big|\langle\varepsilon,\xi\rangle\big|
\frac{1-|y_{1}|^{2}}{|\varepsilon-y_{1}|^{n}}|\varrho_{\varphi_{m}}(y_{1})|dV(y_{1})\\
&\leq&\frac{\|\varphi_{m}\|_{\infty}}{2n^{2}(n+2)}\left[\frac{n+4}{4n(n+2)}\right]^{m-2}|\xi|.
\end{eqnarray*} The proof of this lemma is complete.
\epf

\begin{Lem}{\rm (\cite[Lemma 4]{LX})}\label{LemA}
Let $A$ be an $n \times n$ real $($or complex$)$ matrix with
$|A|\neq0$. Then for any unit vector $\theta\in\partial
\mathbb{B}^{n}$, the inequality
$$|A\theta|\geq\frac{|\det A|}{|A|^{n-1}}
$$
holds.
\end{Lem}

\subsection*{Proof of Theorem \ref{L-B}}
We first prove (a).
By (\ref{eq-ch-3.0}), we have

$$f(x)=P[\varphi_{0}](x)+\sum_{k=1}^{m}(-1)^{k}G_{k}[\varphi_{k}](x),~x\in\mathbb{B}^{n},$$
where $\varphi_{0}=f|_{\mathbb{S}^{n-1}}$. Then

\beqq
D_{f}(x)-D_{f}(0)=D_{P[\varphi_{0}]}(x)-D_{P[\varphi_{0}]}(0)+\sum_{k=1}^{m}(-1)^{k}\left(D_{G_{k}[\varphi_{k}]}(x)-D_{G_{k}[\varphi_{k}]}(0)\right).
\eeqq We split the remaining proof into six steps to complete it.

\bst\label{e1}  The estimate of $|D_{P[\varphi_{0}]}(x)-D_{P[\varphi_{0}]}(0)|$ for $x\in\mathbb{B}^{n}(0,r_{0})$, where $r_{0}=1/2$. \est

  For
$x\in\mathbb{B}^{n}(0,r_{0})$, by Theorem D, we have

\beqq
|D_{P[\varphi_{0}]}(x)-D_{P[\varphi_{0}]}(0)|\leq|D_{P[\varphi_{0}]}(x)|+|D_{P[\varphi_{0}]}(0)|
\leq M_{r_{0}}:=\frac{nM}{2}\frac{2-r_{0}}{1-r_{0}},
\eeqq which, together with Lemma E, implies that

\beq\label{claim-4.1} \nonumber|D_{P[\varphi_{0}]}(x)-D_{P[\varphi_{0}]}(0)|&\leq&
M_{r_{0}}\left(1-\frac{r_{0}^{n-2}(r-|x|)}{(r_{0}+|x|)^{n-1}}\right)\\ \nonumber
&=&M_{r_{0}}\frac{\left({n-1\choose 1}r_{0}^{n-2}+\cdots+{n-1\choose
n-1}|x|^{n-2}+r^{n-2}\right)}{(r_{0}+|x|)^{n-1}}|x|\\ \nonumber
&\leq&M_{r_{0}}\frac{\left((1+r_{0})^{n-1}+r_{0}^{n-2}-r_{0}^{n-1}\right)}{r_{0}^{n-1}}|x|\\
&\leq&\kappa(r_{0})|x|
\eeq where
 $$\kappa(r_{0})=M_{r_{0}}\left(\left(1+\frac{1}{r_{0}}\right)^{n-1}+\frac{1}{r_{0}}-1\right).$$

\bst\label{e3}  The estimate of
$\left|D_{G_{m}[\varphi_{m}]}(x)-D_{G_{m}[\varphi_{m}]}(0)\right|$ for $x\in\mathbb{B}^{n}(0,r_{0})$, where
 $m\geq2$. \est

\bcl\label{set-4.3}
$$\left|D_{G_{m}[\varphi_{m}]}(x)-D_{G_{m}[\varphi_{m}]}(0)\right|\leq\mathscr{L}_{m}(x)$$
and $\lim_{x\rightarrow0}\mathscr{L}_{m}(x)=0$,
 where

$$
\mathscr{L}_{m}(x)=M_{m}\int_{\mathbb{B}^{n}}\left|\nabla_{x}
G(x,y_{1})-\nabla_{x} G(0,y_{1})\right| (1-|y_{1}|^{2})dV(y_{1})
$$ and
$M_{m}=\frac{\|\varphi_{m}\|_{\infty}}{2n}\left[\frac{n+4}{4n(n+2)}\right]^{m-2}$.
\ecl

Now we prove the Claim \ref{set-4.3}. For $\xi\in\mathbb{S}^{n-1}$,
let
$$\mathscr{S}_{\xi}^{m}(x)=\left|(D_{G_{m}[\varphi_{m}]}(x)-D_{G_{m}[\varphi_{m}]}(0))\xi\right|.$$
Then, by (\ref{eq-jj-3}), Lemma F and the Lebesgue Dominated
Convergence Theorem, we have

\beqq
\mathscr{S}_{\xi}^{m}(x)&\leq&\int_{\mathbb{B}^{n}}\left|\langle\nabla_{x}
G(x,y_{1})-\nabla_{x} G(0,y_{1}),\xi\rangle\right|
|\varrho_{\varphi_{m}}(y_{1})|dV(y_{1})\\
&\leq&\int_{\mathbb{B}^{n}}\left|\nabla_{x} G(x,y_{1})-\nabla_{x}
G(0,y_{1})\right||\varrho_{\varphi_{m}}(y_{1})|dV(y_{1})\\
&\leq&\mathscr{L}_{m}(x), \eeqq where $\varrho_{\varphi_{m}}(y_{1})$
is defined in (\ref{eq-jj-2.1}).

On the other hand, by Lemma F and (\ref{eq-hh3}), we see
that
$$\mathscr{L}_{m}(x)\in\mathcal{C}^{0}(\mathbb{B}^{n},\mathbb{R})$$
and

\beqq\lim_{x\rightarrow0}\mathscr{L}_{m}(x)=M_{m}\int_{\mathbb{B}^{n}}\lim_{x\rightarrow0}\left|\nabla_{x}
G(x,y_{1})-\nabla_{x} G(0,y_{1})\right|
(1-|y_{1}|^{2})dV(y_{1})=0.\eeqq The proof of Claim \ref{set-4.3} is
complete.

\bst\label{e2}  For $m\geq2$ and $k\in\{1,\ldots,m-1\}$, we estimate
$\left|D_{G_{k}[\varphi_{k}]}(x)-D_{G_{k}[\varphi_{k}]}(0)\right|$
for $x\in\mathbb{B}^{n}(0,r_{0})$. \est

\bcl\label{set-4.2} For $m\geq2$ and $k\in\{1,\ldots,m-1\}$, we have
$$\left|D_{G_{k}[\varphi_{k}]}(x)-D_{G_{k}[\varphi_{k}]}(0)\right|\leq\mathscr{L}_{k}(x)$$
and $\lim_{x\rightarrow0}\mathscr{L}_{k}(x)=0$,
 where

$$
\mathscr{L}_{k}(x)=\begin{cases} \displaystyle
\|\varphi_{1}\|_{\infty}\int_{\mathbb{B}^{n}}\left|\nabla_{x}
G(x,y_{1})-\nabla_{x} G(0,y_{1})\right| dV(y_{1}),
&\mbox{if}~ k=1,\\
\displaystyle M_{k}\int_{\mathbb{B}^{n}}\left|\nabla_{x}
G(x,y_{1})-\nabla_{x} G(0,y_{1})\right| (1-|y_{1}|^{2})dV(y_{1}),
&\mbox{if}~ 2\leq k\leq m-1,
\end{cases}
$$ and
$M_{k}=\frac{\|\varphi_{k}\|_{\infty}}{2n}\left[\frac{n+4}{4n(n+2)}\right]^{k-2}$.
\ecl In order to prove Claim \ref{set-4.2}, we divide it into two
cases. For $\xi\in\mathbb{S}^{n-1}$, let
$$\mathscr{S}_{\xi}^{k}(x)=\left|\big(D_{G_{k}[\varphi_{k}]}(x)-D_{G_{k}[\varphi_{k}]}(0)\big)\xi\right|.$$

\noindent $\mathbf{Case~ 1.}$ $m\geq3$ and $2\leq k\leq m-1$.



It follows from (\ref{eq-hh2}),  Lemma F and the Lebesgue
Dominated Convergence Theorem that

   \beqq
\mathscr{S}_{\xi}^{k}(x)&\leq&\int_{\mathbb{B}^{n}}\left|\langle\nabla_{x}
G(x,y_{1})-\nabla_{x} G(0,y_{1}),\xi\rangle\right|
|\varrho_{\varphi_{k}}(y_{1})|dV(y_{1})\\
&\leq&\int_{\mathbb{B}^{n}}\left|\nabla_{x} G(x,y_{1})-\nabla_{x}
G(0,y_{1})\right||\varrho_{\varphi_{k}}(y_{1})|dV(y_{1})\\
&\leq&\mathscr{L}_{k}(x), \eeqq where $\varrho_{\varphi_{k}}(y_{1})$
is defined in (\ref{eq-hh-1}).

Next we prove $\lim_{x\rightarrow0}\mathscr{L}_{k}(x)=0$.



By Lemma F and (\ref{eq-hh3}), we see that
$$\mathscr{L}_{k}(x)\in\mathcal{C}^{0}(\mathbb{B}^{n},\mathbb{R})$$
and

\beqq\lim_{x\rightarrow0}\mathscr{L}_{k}(x)=M_{k}\int_{\mathbb{B}^{n}}\lim_{x\rightarrow0^{+}}\left|\nabla_{x}
G(x,y_{1})-\nabla_{x} G(0,y_{1})\right|
(1-|y_{1}|^{2})dV(y_{1})=0.\eeqq

\noindent $\mathbf{Case~ 2.}$ $m=2$ and $k=1$.

In this case, we have
 \beqq
\mathscr{S}_{\xi}^{1}(x)\leq\|\varphi_{1}\|_{\infty}\int_{\mathbb{B}^{n}}\left|\langle\nabla_{x}
G(x,y_{1})-\nabla_{x} G(0,y_{1}),\xi\rangle\right| dV(y_{1})
 \leq\mathscr{L}_{1}(x). \eeqq
It follows from Lemma F and (\ref{eq-hh3})  that
$$\mathscr{L}_{1}(x)\in\mathcal{C}^{0}(\mathbb{B}^{n},\mathbb{R})$$
and

\beqq\lim_{x\rightarrow0}\mathscr{L}_{1}(x)=\|\varphi_{1}\|_{\infty}\int_{\mathbb{B}^{n}}\lim_{x\rightarrow0}\left|\nabla_{x}
G(x,y_{1})-\nabla_{x} G(0,y_{1})\right| dV(y_{1})=0.\eeqq The proof
of Claim \ref{set-4.2} is finished.

\bst\label{e4}  The estimate of $|D_{f}(0)\xi|$, where
$\xi\in\mathbb{S}^{n-1}$. \est

By Theorem D, Lemmas  \ref{eq-jj-1} and \ref{eq-jj-2}, we have

\be\label{ll-1} |D_{f}(0)|\leq\left|D_{P[\varphi_{0}]}(0)\right|+\sum_{k=1}^{m}\left|D_{G_{k}[\varphi_{k}]}(0)\right|\leq M^{\ast},
\ee where $M^{\ast}=\frac{n}{2}M+\frac{n}{n+1}\|\varphi_{1}\|_{\infty}+\sum_{k=2}^{m}\frac{\|\varphi_{k}\|_{\infty}}{2n}\left[\frac{n+4}{4n(n+2)}\right]^{k-2}\delta(n)$
and $\delta(n)$ is defined in Lemma {\rm \ref{eq-jj-1}}.

For any $\xi\in\mathbb{S}^{n-1}$, it follows from (\ref{ll-1}) and
Lemma G that

\be\label{ll-2}
|D_{f}(0)\xi|\geq\frac{J_{f}(0)}{|D_{f}(0)|^{n-1}}\geq
\frac{1}{(M^{\ast})^{n-1}}. \ee 

\bst\label{e5} We will show that there is a constant
$r_{1}\in(0,r_{0})$ such that $f$ is   injective in
$\mathbb{B}^{n}(0,r_{1})$. \est

In order to prove the injection of $f$ in $\mathbb{B}^{n}(0,r_{1})$,
let
$$\mathscr{G}(|x|)=\frac{1}{(M^{\ast})^{n-1}}-\kappa(r_{0})|x|-
\max_{0\leq|\zeta|\leq|x|}\left(\sum_{k=1}^{m}\mathscr{L}_{k}(\zeta)\right),$$
$x\in\mathbb{B}^{n}(0,r_{0})$. Since $\mathscr{G}(|x|)$ is
monotonous and continuous in $\mathbb{B}^{n},$
$$\lim_{|x|\rightarrow0^{+}}\mathscr{G}(|x|)=\frac{1}{(M^{\ast})^{n-1}}>0~\mbox{and}~\lim_{|x|\rightarrow r_{0}}\mathscr{G}(|x|)<0,$$ we see that
there is an unique constant $r_{1}\in(0,r_{0}]$ such that
$\mathscr{G}(r_{1})=0$. For any $x',x''\in\mathbb{B}^{n}(0,r_{1})$,
we use $[x',x'']$ to denote the segment from $x'$ to $x''$ with the
endpoints $x'$ and $x''$. Hence by (\ref{claim-4.1}), Claims
\ref{set-4.3}, \ref{set-4.2}  and  (\ref{ll-2}),  we have

\beqq
|f(x')-f(x'')|&\geq&\left|\int_{[x',x'']}D_{f}(0)dx\right|-\int_{[x',x'']}\left|D_{f}(x)-D_{f}(0)\right||dx|\\ \nonumber
&\geq&\left|\int_{[x',x'']}D_{f}(0)dx\right|-\int_{[x',x'']}\left|D_{P[\varphi_{0}]}(x)-D_{P[\varphi_{0}]}(0)\right||dx|\\ \nonumber
&&-\int_{[x',x'']}\sum_{k=1}^{m}\left|D_{G_{k}[\varphi_{k}]}(x)-D_{G_{k}[\varphi_{k}]}(0)\right||dx|\\ \nonumber
&\geq&\int_{[x',x'']}\mathscr{G}(x)|dx|>\int_{[x',x'']}\mathscr{G}(r_{1})|dx|\\ \nonumber
&=&|x'-x''|\Bigg(\frac{1}{(M^{\ast})^{n-1}}-\kappa(r_{0})r_{1}
-\max_{0\leq|\zeta|\leq
r_{1}}\left(\sum_{k=1}^{m}\mathscr{L}_{k}(\zeta)\right)\Bigg)\\
&=&0,
\eeqq which yields that $f(x')\neq f(x'')$. Thus, from the
arbitrariness of $x'$ and $x''$, the injection of $f$ follows.

\bst\label{e6} We will prove that the image
$f(\mathbb{B}^{n}(0,r_{1}))$  contains a  ball
$\mathbb{B}^{n}(0,R_{1})$, where
$$R_{1}\geq\frac{\kappa(r_{0})r_{1}}{2}=\frac{3n(1+3^{n-1})}{4}r_{1}M.$$
\est

 To reach this
goal, let $\varsigma\in\partial\mathbb{B}^{n}(0,r_{1})$. Then we
infer from (\ref{claim-4.1}), Claims  \ref{set-4.3},  \ref{set-4.2} and
(\ref{ll-2}) that

\beqq
|f(\varsigma)-f(0)|&\geq&\left|\int_{[0,\varsigma]}D_{f}(0)dx\right|-\int_{[0,\varsigma]}\left|D_{f}(x)-D_{f}(0)\right||dx|\\
&\geq&\left|\int_{[0,\varsigma]}D_{f}(0)dx\right|-\int_{[0,\varsigma]}\left|D_{P[\varphi_{0}]}(x)-D_{P[\varphi_{0}]}(0)\right||dx|\\
&&-\int_{[0,\varsigma]}\sum_{k=1}^{m}\left|D_{G_{k}[\varphi_{k}]}(x)-D_{G_{k}[\varphi_{k}]}(0)\right||dx|
\geq\int_{[0,\varsigma]}\mathscr{G}(x)|dx|\\
&=&r_{1}\Bigg(\frac{1}{(M^{\ast})^{n-1}}-\frac{\kappa(r_{0})}{2}r_{1}-\max_{0\leq|\zeta|\leq
r_{1}}\left(\sum_{k=1}^{m}\mathscr{L}_{k}(\zeta)\right)\Bigg)\\
&=&\frac{\kappa(r_{0})}{2}r_{1},~\mbox{(by~ $\mathscr{G}(r_{1})=0$)},
\eeqq
which implies that $f(\mathbb{B}^{n}(0,r_{1}))$ contains a
 ball $\mathbb{B}^{n}(0,R_{1})$ with $$R_{1}\geq\frac{\kappa(r_{0})r_{1}}{2}=\frac{3n(1+3^{n-1})}{4}r_{1}M.$$

 Next we prove (b). If $m=1$, then, by (\ref{po-1}), we have
 $$f(x)=P[\varphi_{0}](x)-G_{1}[\varphi_{1}](x),~x\in\mathbb{B}^{n},$$
where $\varphi_{0}=f|_{\mathbb{S}^{n-1}}$ and $G_{1}[\varphi_{1}]$ is defined in (\ref{po-2}).
In the following, we will use the similar reasoning as in the proof of (a) to prove (b).
Let $r_{1}$ be an unique solution of the following   equation: \beqq
&&\frac{1}{\left(\frac{nM}{2}+\frac{n\|\varphi_{1}\|_{\infty}}{n+1}\right)^{n-1}}-
\kappa(r_{0})r_{1}\\
&-&\|\varphi_{1}\|_{\infty} \max_{0\leq|x|\leq
r_{1}}\left(\int_{\mathbb{B}^{n}}\left|\nabla_{x}
G(x,y_{1})-\nabla_{x} G(0,y_{1})\right| dV(y_{1})\right)=0.\eeqq
 Then, for any $x',x''\in\mathbb{B}^{n}(0,r_{1})$, we have

  \beqq
|f(x')-f(x'')|&\geq&\left|\int_{[x',x'']}D_{f}(0)dx\right|-\int_{[x',x'']}\left|D_{f}(x)-D_{f}(0)\right||dx|\\ \nonumber
&\geq&\left|\int_{[x',x'']}D_{f}(0)dx\right|-\int_{[x',x'']}\left|D_{P[\varphi_{0}]}(x)-D_{P[\varphi_{0}]}(0)\right||dx|\\ \nonumber
&&-\int_{[x',x'']}\left|D_{G_{1}[\varphi_{1}]}(x)-D_{G_{1}[\varphi_{1}]}(0)\right||dx|\\ \nonumber
&>&|x'-x''|\Bigg(\frac{1}{\left(\frac{nM}{2}+\frac{n\|\varphi_{1}\|_{\infty}}{n+1}\right)^{n-1}}
-\kappa(r_{0})r_{1}\\ \nonumber
&&-\|\varphi_{1}\|_{\infty} \max_{0\leq|x|\leq
r_{1}}\left(\int_{\mathbb{B}^{n}}\left|\nabla_{x}
G(x,y_{1})-\nabla_{x} G(0,y_{1})\right| dV(y_{1})\right)\\ \nonumber
&=&0,
\eeqq which yields that $f$ is injective in $\mathbb{B}^{n}(0,r_{1})$.
 Therefore, for any $\varsigma\in\partial\mathbb{B}^{n}(0,r_{1})$, we see that

  \beqq
|f(\varsigma)-f(0)|&\geq&\left|\int_{[0,\varsigma]}D_{f}(0)dx\right|-\int_{[0,\varsigma]}\left|D_{f}(x)-D_{f}(0)\right||dx|\\
&\geq&\left|\int_{[0,\varsigma]}D_{f}(0)dx\right|-\int_{[0,\varsigma]}\left|D_{P[\varphi_{0}]}(x)-D_{P[\varphi_{0}]}(0)\right||dx|\\
&&-\int_{[0,\varsigma]}\left|D_{G_{1}[\varphi_{1}]}(x)-D_{G_{1}[\varphi_{1}]}(0)\right||dx|
\\
&\geq&\frac{\kappa(r_{0})}{2}r_{1}.
\eeqq
  The proof of this theorem is complete. \qed

\section{The Lipschitz  continuity of quasiconformal self-mappings  satisfying  polyharmonic equations}\label{csw-sec4}


\begin{lem}\label{lem-mck}
Suppose that $n\geq3$, $m\in\mathbb{N}\setminus\{1\}$,
$\varphi_{m}\in\mathcal{C}(\overline{\mathbb{B}^{n}},\mathbb{R}^{n})$
and $\varphi_{k}\in \mathcal{C}(\mathbb{S}^{n-1},\mathbb{R}^{n})$
for $k\in\{0,1,\ldots,m-1\}$. Let $f$ be a mapping of
$\overline{\mathbb{B}^{n}}$ onto itself satisfying {\rm (\ref{eq-ch-1})}
and the boundary conditions
$\Delta^{m-1}f=\varphi_{m-1},~\ldots,~\Delta^{1}f=\varphi_{1}$ on
$\mathbb{S}^{n-1}$. In addition, let $f$ be Lipschitz continuous in
$\mathbb{B}^{n}$ satisfying $|f(x)|\rightarrow1$ as
$|x|\rightarrow1^{-}$, where $x\in\mathbb{B}^{n}$. Then, for almost
every $t\in\mathbb{S}^{n-1}$, the following limits exist:

\be\label{jj-5}
D_{f}(t):=\lim_{r\rightarrow1^{-}}D_{f}(rt)~\mbox{and}~J_{f}(t):=\lim_{r\rightarrow1^{-}}J_{f}(rt).\ee

Further, for $\varphi_{0}:=f|_{\mathbb{S}^{n-1}}$ and
$x(\theta)=\varphi_{0}(T(\theta)):=\varphi_{0}(t)$, we have
\beq\label{M-1.1}J_{f}(t)&\leq&\frac{M_{x}(\theta)}{M_{T}(\theta)}\bigg\{\int_{\mathbb{S}^{n-1}}\frac{|\varphi_{0}(t)-\varphi_{0}(\zeta)|^{2}}{|\zeta-t|^{n}}d\sigma(\zeta)
+\frac{\|\varphi_{1}\|_{\infty}}{n}\\ \nonumber
&&+\sum_{k=2}^{m}\frac{\|\varphi_{k}\|_{\infty}}{n^{2}(n+2)}\left[\frac{n+4}{4n(n+2)}\right]^{k-2}\bigg\}
\eeq and
\beq\label{M-1.2}J_{f}(t)&\geq&\frac{M_{x}(\theta)}{M_{T}(\theta)}\bigg\{\int_{\mathbb{S}^{n-1}}\frac{|\varphi_{0}(t)-\varphi_{0}(\zeta)|^{2}}{|\zeta-t|^{n}}d\sigma(\zeta)
-\frac{\|\varphi_{1}\|_{\infty}}{n}\\ \nonumber
&&-\sum_{k=2}^{m}\frac{\|\varphi_{k}\|_{\infty}}{n^{2}(n+2)}\left[\frac{n+4}{4n(n+2)}\right]^{k-2}\bigg\},
\eeq where  $M_{x}(\theta)$ and $M_{T}(\theta)$ are the square roots
of Gram determinants of $D_{x}$ and $D_{T}$, respectively.
\end{lem}

Before the proof of Lemma \ref{lem-mck}, let us recall the following
result.

\begin{Lem}\label{K-1}{\rm (\cite[Lemma 2.1]{K-2011})}
Let $u=P[f]$ be a harmonic mapping of $\mathbb{B}^{n}$ into
$\mathbb{R}^{n}$, and assume that its derivative $v=D_{u}$ is
bounded in $\mathbb{B}^{n}$ (or equivalently, let $u$ be Lipschitz
continuous), where $u|_{\mathbb{S}^{n-1}}=f\in
L^{1}(\mathbb{S}^{n-1})$. Then there exists a mapping $A\in
L^{\infty}(\mathbb{S}^{n-1})$ defined in the $\mathbb{S}^{n-1}$ such
 that $D_{f}(x)=P[A](x)$ and for almost every
 $\eta\in\mathbb{S}^{n-1}$ there holds the relation
 $$\lim_{r\rightarrow1^{-}}D_{f}(r\eta)=A(\eta).$$  Moveover,  the
 function $f\circ T$ is differentiable almost everywhere in
 $Q^{n-1}$ and there holds $$A(T(\theta))D_{f}(\theta)=D_{f\circ
 T}(\theta),$$ where
 $\theta=(\theta_{1},\theta_{2},\ldots,\theta_{n-1})$, $Q^{n-1}$ and $T$
 are
 defined in the part of  {\rm \ref{sbcsw-sec2.3}}.
\end{Lem}

\subsection*{The proof of Lemma \ref{lem-mck}}
 We first prove the existence of the two limits in (\ref{jj-5}).
By Lemmas \ref{eq-jj-1} and  \ref{eq-jj-2}, we see that for any
 $t\in\mathbb{S}^{n-1}$,
\be\label{jj-5z}
\lim_{r\rightarrow1^{-}}D_{G_{k}[\varphi_{k}]}(rt)=D_{G_{k}[\varphi_{k}]}(t)\ee
and $G_{k}[\varphi_{k}]$ are Lipschitz continuous in
$\mathbb{B}^{n}$, where $k\in\{1,\ldots,m\}$. Since $f$ is Lipschitz
continuous in $\mathbb{B}^{n}$, we see that
$$P[\varphi_{0}](x)=f(x)-\sum_{k=1}^{m}(-1)^{k}G_{k}[\varphi_{k}](x)$$
are also Lipschitz continuous in $\mathbb{B}^{n}$, where
$\varphi_{0}=f|_{\mathbb{S}^{n-1}}$. It follows from Lemma H
that, for almost every $t\in\mathbb{S}^{n-1}$,
$$D_{P[\varphi_{0}]}(t)=\lim_{r\rightarrow1^{-}}D_{P[\varphi_{0}]}(rt)$$
does exist, which, together with (\ref{jj-5z}), guarantees that for
almost every $t\in\mathbb{S}^{n-1}$,
\be\label{jj-6}D_{f}(t)=\lim_{r\rightarrow1^{-}}D_{f}(rt)\ee also
exists.

By (\ref{jj-6}) and $J_{f}=\det D_{f},$ we conclude that
$$J_{f}(t)=\lim_{r\rightarrow1^{-}}J_{f}(rt)$$ exists for almost every
$t\in\mathbb{S}^{n-1}$.

Next we estimate $J_{f}(t)$. It follows from (\ref{jj-6}) that the
mapping $x$, $x(\theta)=\varphi_{0}(T(\theta))$, defines the outer
normal vector field $\mathbf{n}_{x}$ almost everywhere in
$\mathbb{S}^{n-1}$ at the point
$x(\theta)=\varphi_{0}(T(\theta))=(x_{1},\ldots ,x_{n})'$ by the
formula

$$\mathbf{n}_{x}(x(\theta))=x_{\theta_{1}}\times\cdots\times x_{\theta_{n-2}}\times x_{\theta_{n-1}}=\left(\begin{array}{cccc}
e_{1}\;\;\;\;\; \; e_{2}\;\;\;\;\;\;\cdots\;\;\;\;\;\;
  e_{n}\\
 x_{1\theta_{1}}\;\;\;\; x_{2\theta_{1}}\;\;\;\;\cdots\;\;\;\;
 x_{n\theta_{1}}\\
\vdots \\
 x_{1\theta_{n-2}}\;\; x_{2\theta_{n-2}}\;\;
\cdots\;\;
x_{n\theta_{n-2}}\\
 x_{1\theta_{n-1}}\;\; x_{2\theta_{n-1}}\;\;\cdots\;\;
  x_{n\theta_{n-1}}
\end{array}\right),
$$
where $e_{1}=(1,0,\ldots,0)'$, $\ldots$, $e_{n}=(0,0,\ldots,1)'$ and
$T$ is defined in the part of  \ref{sbcsw-sec2.3}. Let
$f(S(r,\theta))=y=(y_{1},\ldots,y_{n})'$, where $S$ is defined in
\ref{sbcsw-sec2.3}.

By (\ref{jj-6}), for $i\in\{1,\ldots,n\}$ and
$j\in\{1,\ldots,n-1\}$, we have

$$\lim_{r\rightarrow1^{-}}y_{i\theta_{j}}(r,\theta)=x_{i\theta_{j}}(\theta)$$
and
$$\lim_{r\rightarrow1^{-}}y_{ir}(r,\theta)=
\lim_{r\rightarrow1^{-}}\frac{x_{i}(\theta)-y_{i}(r,\theta)}{1-r},$$
which imply that

\beq\label{eq-jj-7}\nonumber \lim_{r\rightarrow1^{-}}J_{f\circ
S}(r,\theta)&=&
\lim_{r\rightarrow1^{-}}\left\langle\frac{x-y}{1-r},x_{\theta_{1}}\times\cdots\times
 x_{\theta_{n-1}}\right\rangle\\ \nonumber
&=&\lim_{r\rightarrow1^{-}}\left\langle\frac{x-P[\varphi_{0}]}{1-r},x_{\theta_{1}}\times\cdots\times
x_{\theta_{n-1}}\right\rangle-\sum_{k=1}^{m}(-1)^{k}\mathcal{X}_{k}\\
\nonumber
 &=&\lim_{r\rightarrow1^{-}}\int_{\mathbb{S}^{n-1}}\frac{1+r}{|\zeta-rt|^{n}}\left\langle x-\varphi_{0}(\zeta),x_{\theta_{1}}\times\cdots\times
 x_{\theta_{n-1}}\right\rangle\;
d\sigma(\zeta)\\ \nonumber &&-\sum_{k=1}^{m}(-1)^{k}\mathcal{X}_{k}\\
 &=&\lim_{r\rightarrow1^{-}}\int_{\mathbb{S}^{n-1}}\frac{1+r}{|\zeta-S(r,\theta)|^{n}}\left\langle \varphi_{0}(T(\theta))-\varphi_{0}(\zeta),\mathbf{n}_{\varphi_{0}\circ T}(T(\theta))\right\rangle\;
d\sigma(\zeta)\\ \nonumber &&-\sum_{k=1}^{m}(-1)^{k}\mathcal{X}_{k},
\eeq where

\be\label{eq-jj-7.1}\mathcal{X}_{k}=\lim_{r\rightarrow1^{-}}\left\langle\frac{G_{k}[\varphi_{k}](x)}{1-r},x_{\theta_{1}}\times\cdots\times
 x_{\theta_{n-1}}\right\rangle.\ee

Since
$$\mathbf{n}_{x}(x(\theta))=M_{x}(\theta)\cdot\varphi_{0}(T(\theta)),$$
by (\ref{eq-jj-7}), we see that

\beq\label{eq-jj-8}\nonumber \lim_{r\rightarrow1^{-}}J_{f\circ
S}(r,\theta)&=&
M_{x}(\theta)\lim_{r\rightarrow1^{-}}\int_{\mathbb{S}^{n-1}}\frac{1+r}{|\zeta-S(r,\theta)|^{n}}\left\langle
\varphi_{0}(T(\theta))-\varphi_{0}(\zeta),\varphi_{0}(T(\theta))\right\rangle
d\sigma(\zeta)\\ \nonumber
&&-\sum_{k=1}^{m}(-1)^{k}\mathcal{X}_{k}\\ &=& M_{x}(\theta)
\lim_{r\rightarrow1^{-}}\int_{\mathbb{S}^{n-1}}\frac{\left|\varphi_{0}(T(\theta))-\varphi_{0}(\zeta)\right|^{2}}{|\zeta-S(r,\theta)|^{n}}d\sigma(\zeta)
-\sum_{k=1}^{m}(-1)^{k}\mathcal{X}_{k}.
 \eeq
In the following, we will demonstrate the estimate of
$|\mathcal{X}_{k}|$ for $k\in\{1,\ldots,m\}$.

 \noindent
$\mathbf{Case~ 1.}$   $k=1$. Then, by  (\ref{eq-jj-7.1}), we get

\begin{eqnarray*}
\mathcal{X}_{1}&=&\lim_{x\rightarrow
t\in\mathbb{S}^{n-1}}\left\langle\frac{G_{1}[\varphi_{1}](x)}{1-|x|},x_{\theta_{1}}\times\cdots\times
 x_{\theta_{n-1}}\right\rangle\\ \nonumber&=&
M_{x}(\theta)\int_{\mathbb{B}^{n}}\lim_{x\rightarrow
t}\frac{G(x,y_{1})}{1-|x|}\left\langle\int_{\mathbb{S}^{n-1}}P(y_{1},\zeta)\varphi_{1}(\zeta)d\sigma(\zeta),\varphi_{0}(t)\right\rangle dV(y_{1})\\
\nonumber
&=&\frac{M_{x}(\theta)}{\omega_{n-1}}\int_{\mathbb{B}^{n}}P(y_{1},t)\left\langle\int_{\mathbb{S}^{n-1}}P(y_{1},\zeta)
\varphi_{1}(\zeta)d\sigma(\zeta),\varphi_{0}(t)\right\rangle
dV(y_{1}),
\end{eqnarray*}  which implies that

\beq\label{eq-jj9}
|\mathcal{X}_{1}|&\leq&\frac{M_{x}(\theta)}{\omega_{n-1}}\|\varphi_{1}\|_{\infty}\int_{\mathbb{B}^{n}}P(y_{1},t)dV(y_{1})\\
\nonumber
&=&M_{x}(\theta)\|\varphi_{1}\|_{\infty}\int_{0}^{1}\left(\rho^{n-1}\int_{\mathbb{S}^{n-1}}P(\rho\zeta,t)d\sigma(\zeta)\right)d\rho\\
\nonumber &=&\frac{M_{x}(\theta)\|\varphi_{1}\|_{\infty}}{n}. \eeq

\noindent $\mathbf{Case~ 2.}$   $k\in\{2,\ldots,m-1\}$ and $m\geq3$.
In this case, by (\ref{eq-jj-7.1}),  we have

\begin{eqnarray*}
 \mathcal{X}_{k}&=&\lim_{x\rightarrow
t\in\mathbb{S}^{n-1}}\left\langle\frac{G_{k}[\varphi_{k}](x)}{1-|x|},x_{\theta_{1}}\times\cdots\times
 x_{\theta_{n-1}}\right\rangle\\ \nonumber&=&
M_{x}(\theta)\int_{\mathbb{B}^{n}}\lim_{x\rightarrow
t}\frac{G(x,y_{1})}{1-|x|}\left\langle\varrho_{\varphi_{k}}(y_{1}),\varphi_{0}(t)\right\rangle dV(y_{1})\\
\nonumber&=&\frac{M_{x}(\theta)}{\omega_{n-1}}\int_{\mathbb{B}^{n}}P(y_{1},t)\left\langle\varrho_{\varphi_{k}}(y_{1}),\varphi_{0}(t)\right\rangle
dV(y_{1}),
\end{eqnarray*} which, together with (\ref{eq-hh2}),
gives that

\beq\label{eq-jj-11}
|\mathcal{X}_{k}|&\leq&\frac{M_{x}(\theta)}{\omega_{n-1}}\int_{\mathbb{B}^{n}}P(y_{1},t)
|\varrho_{\varphi_{k}}(y_{1})| dV(y_{1})\\ \nonumber
&\leq&\frac{M_{x}(\theta)\|\varphi_{k}\|_{\infty}}{2n\omega_{n-1}}\left[\frac{n+4}{4n(n+2)}\right]^{k-2}\int_{\mathbb{B}^{n}}P(y_{1},t)(1-|y_{1}|^{2})
 dV(y_{1})\\ \nonumber
 &=&\frac{M_{x}(\theta)\|\varphi_{k}\|_{\infty}}{n^{2}(n+2)}\left[\frac{n+4}{4n(n+2)}\right]^{k-2},\eeq
where $\varrho_{\varphi_{k}}(y_{1})$ is defined in (\ref{eq-hh-1}).

\noindent $\mathbf{Case~ 3.}$ $k=m$. Then it follows from
(\ref{eq-jj-3}) and (\ref{eq-jj-7.1}) that

\beq\label{eq-jj-12}
 |\mathcal{X}_{m}|&\leq&\frac{M_{x}(\theta)}{\omega_{n-1}}\left|\int_{\mathbb{B}^{n}}P(y_{1},t)
 \left\langle\varrho_{\varphi_{m}}(y_{1}),\varphi_{0}(t)\right\rangle
dV(y_{1})\right|\\ \nonumber
&\leq&\frac{M_{x}(\theta)}{\omega_{n-1}}\int_{\mathbb{B}^{n}}P(y_{1},t)
|\varrho_{\varphi_{m}}(y_{1})| dV(y_{1})\\ \nonumber
&\leq&\frac{M_{x}(\theta)\|\varphi_{m}\|_{\infty}}{2n\omega_{n-1}}\left[\frac{n+4}{4n(n+2)}\right]^{m-2}\int_{\mathbb{B}^{n}}P(y_{1},t)(1-|y_{1}|^{2})
 dV(y_{1})\\ \nonumber
 &=&\frac{M_{x}(\theta)\|\varphi_{m}\|_{\infty}}{n^{2}(n+2)}\left[\frac{n+4}{4n(n+2)}\right]^{m-2},
\eeq where $\varrho_{\varphi_{m}}(y_{1})$ is defined in
(\ref{eq-jj-2.1}). Hence (\ref{M-1.1}) and (\ref{M-1.2}) follow from
(\ref{eq-jj-8}), (\ref{eq-jj9}), (\ref{eq-jj-11}), (\ref{eq-jj-12})
and $$J_{f\circ
S}(r,\theta)=r^{n-1}J_{f}(rT(\theta))M_{T}(\theta).$$ The proof of
this lemma is complete.
 \qed

\begin{Lem}\label{K-2.2}{\rm (\cite[Lemma 2.2]{K-2011})}
Let $u$ be a harmonic Lipschitz continuous mapping defined in
$\mathbb{B}^{n}$. Suppose that $D_{u}$ exists almost everywhere in
$\mathbb{S}^{n-1}$. Then for $x\in\mathbb{B}^{n}$,
$$|D_{u}(x)|\leq {\rm ess
~sup}_{|\eta|=1}|D_{u}(\eta)|.$$
\end{Lem}

\begin{Lem}\label{K-2.8}{\rm (\cite[Lemma 4.8]{K3})}
Let $\mathcal{A}:~\mathbb{R}^{n}\rightarrow\mathbb{R}^{n}$ be a
linear operator such that $\mathcal{A}=[a_{ij}]_{i,j=1,\ldots,n}$.
If $\mathcal{A}$ is $K$-quasiconformal, then  the following sharp
inequalities hold:

$$K^{1-n}|\mathcal{A}|^{n-1}\left|x_{1}\times\cdots\times x_{n-1}\right|\leq\left|\mathcal{A}x_{1}\times\cdots\times \mathcal{A}x_{n-1}\right|
\leq|\mathcal{A}|^{n-1}\left|x_{1}\times\cdots\times
x_{n-1}\right|.$$

\end{Lem}

\begin{Lem}\label{K-3.7}{\rm (\cite[Corollary 3.7]{K3})}
Assume that $u:~\overline{\mathbb{B}^{n}}\rightarrow\mathbb{R}^{n}$
is  a  $K$-quasiregular, twice differentiable mapping, continuous on
$\overline{\mathbb{B}^{n}}$, and that
$u|_{\mathbb{S}^{n-1}}\in\mathcal{C}^{1,\alpha}$. If, in addition,
$u$ satisfies the differential inequality

$$|\Delta u|\leq a|D_{u}|^{2}+b$$ for some positive constants $a$
and $b$, then $|D_{u}|$ is bounded and $u$ is Lipschitz continuous.
\end{Lem}

The following is the so-called Mori's Theorem of quasiconformal
mappings defined in $\mathbb{B}^{n}$ (see \cite{FV}).

\begin{Thm}\label{FV-Th-H}
If $u$ is a $K$-quasiconformal self-mapping of $\mathbb{B}^{n}$ with
$u(0)=0$, then there exists a constant $q(n,K)$, satisfying the
condition $q(n,K)\rightarrow1$ as $K\rightarrow1$, such that, for
any $x,y\in\mathbb{B}^{n}$, $$|u(x)-u(y)|\leq
q(n,K)|x-y|^{K^{1/(1-n)}}.$$ Moreover, the mapping
$u(x)=|x|^{-1+K^{1/(1-n)}}x$ shows that the exponent $K^{1/(1-n)}$
is optimal in the class of arbitrary $K$-quasiconformal
homeomorphism from $\mathbb{B}^{n}$  onto itself.
\end{Thm}

\subsection*{The proof of Theorem \ref{thm-2.1}}
 Let's begin  the proof of this theorem with
the following claim.

\bcl\label{Ca-1} The limits $$\lim_{x\rightarrow
\xi\in\mathbb{S}^{n-1},x\in\mathbb{B}^{n}}D_{f}(x)\;\;
\mbox{and}\;\; \lim_{x\rightarrow
\xi\in\mathbb{S}^{n-1},x\in\mathbb{B}^{n}}J_{f}(x)$$ exist almost
everywhere in $\mathbb{S}^{n-1}$. \ecl

In order to prove the existence of these two limits, we need to
obtain the upper bound of $|\Delta f(x)|$ in $\mathbb{B}^{n}$, and
we divide it into two cases to estimate.

\noindent $\mathbf{Case~ 1.}$  $m=2.$

By  \cite[pp. 118-120]{Ho} (see also \cite{K3,K4}), we have that for
$x\in\mathbb{B}^{n}$, \beqq\label{eq-yy2}\Delta
f(x)=P[\varphi_{1}](x)-\int_{\mathbb{B}^{n}}G(x,\zeta)\varphi_{2}(\zeta)dV(\zeta).
\eeqq It follows from Lemma B  that \be\label{eq-x1.1}
|\Delta f(x)|\leq
|P[\varphi_{1}](x)|+\|\varphi_{2}\|_{\infty}\int_{\mathbb{B}^{n}}|G(x,\zeta)|dV(\zeta)
\leq
\|\varphi_{1}\|_{\infty}+\frac{\|\varphi_{2}\|_{\infty}}{2n}<\infty.\ee

\noindent $\mathbf{Case~ 2.}$  $m\geq3.$

By (\ref{eq-ch-3.0}), we have that for $x\in \mathbb{B}^{n}$,
\begin{eqnarray*}
\Delta
f(x)=P[\varphi_{1}](x)+\sum_{k=1}^{m-1}(-1)^{k}G_{k}[\varphi_{k+1}](x),
\end{eqnarray*} where
\begin{eqnarray*}
G_{k}[\varphi_{k+1}](x)&=&
\int_{\mathbb{B}^{n}}\cdots\int_{\mathbb{B}^{n}}G(x,\xi_{1})\cdots
G(\xi_{k-1},\xi_{k})\\ \nonumber
&&\times\left(\int_{\mathbb{S}^{n-1}}P(\xi_{k},\xi)\varphi_{k+1}(\xi)d\sigma(\xi)\right)dV(\xi_{k})\cdots
dV(\xi_{1})
\end{eqnarray*}for $k\in\{1,\ldots,m-2\}$, and

\begin{eqnarray*}
 G_{m-1}[\varphi_{m}](x)&=&
\int_{\mathbb{B}^{n}}\cdots\int_{\mathbb{B}^{n}}G(x,\zeta_{1})\cdots G(\zeta_{m-3},\zeta_{m-2})\\
\nonumber
&&\times\left(\int_{\mathbb{B}^{n}}G(\zeta_{m-2},\zeta_{m-1})\varphi_{m}(\zeta_{m-1})dV(\zeta_{m-1})\right)dV(\zeta_{m-2})\cdots
dV(\zeta_{1}).
\end{eqnarray*}

 For $x\in\mathbb{B}^{n}$ and
$k\in\{1,\ldots,m-2\}$, by  (\ref{eq-ch-4g}), we obtain

\begin{eqnarray*}
|G_{k}[\varphi_{k+1}](x)|&\leq&\|\varphi_{k+1}\|_{\infty}
\int_{\mathbb{B}^{n}}\cdots\int_{\mathbb{B}^{n}}|G(x,\xi_{1})|\cdots
|G(\xi_{k-1},\xi_{k})|dV(\xi_{k})\cdots dV(\xi_{1})\\
&\leq&\frac{\|\varphi_{k+1}\|_{\infty}}{2n}\left(\frac{n+4}{4n(n+2)}\right)^{k-1}(1-|x|^{2})\\
&\leq&\frac{\|\varphi_{k+1}\|_{\infty}}{2n}\left(\frac{n+4}{4n(n+2)}\right)^{k-1},
\end{eqnarray*}
and, by (\ref{eq-ch-6g}), we have

\begin{eqnarray*}
 |G_{m-1}[\varphi_{m}](x)|&=&\|\varphi_{m}\|_{\infty}
\int_{\mathbb{B}^{n}}\cdots\int_{\mathbb{B}^{n}}|G(x,\zeta_{1})|\cdots
|G(\zeta_{m-3},\zeta_{m-2})|dV(\zeta_{m-2})\cdots
dV(\zeta_{1})\\
&\leq&\frac{\|\varphi_{m}\|_{\infty}}{2n}\left(\frac{n+4}{4n(n+2)}\right)^{m-2}(1-|x|^{2})\\
&\leq&\frac{\|\varphi_{m}\|_{\infty}}{2n}\left(\frac{n+4}{4n(n+2)}\right)^{m-2},
\end{eqnarray*}
which imply that

\beq\label{eq-x1.01} |\Delta
f(x)|&=&|P[\varphi_{1}](x)|+\sum_{k=1}^{m-1}|G_{k}[\varphi_{k+1}](x)|\\
\nonumber
&\leq&\|\varphi_{1}\|_{\infty}+\sum_{k=1}^{m-1}\frac{\|\varphi_{k+1}\|_{\infty}}{2n}\left(\frac{n+4}{4n(n+2)}\right)^{k-1}<\infty.
\eeq

Since $f$ is a $K$-quasiconformal self-mapping of $\mathbb{B}^{n}$,
we see that $f$ can be extended to the homeomorphism of
$\overline{\mathbb{B}^{n}}$ onto itself. Hence Claim \ref{Ca-1}
follows from (\ref{eq-x1.01}), Lemmas K and \ref{lem-mck}.

In the following, for convenience, let
$$ C_{2}(K,\varphi_{1},\cdots,\varphi_{n})=\sup_{x\in\mathbb{B}^{n}}|D_{f}(x)|.$$

Since for almost all $x_1$ and $x_2\in \mathbb{B}^{n}$,
\beqq\label{sun-7} |f(x_1)-f(x_2)=\Big|\int_{[x_1,x_2]}D_{f}(x)dx
\Big|\leq C_{2}(K,\varphi_{1},\cdots,\varphi_{n}) |x_1-x_2|, \eeqq
we see that, to prove the Lipschitz continuity of $f$, it suffices
to estimate the quantity $C_{2}(K,\varphi_{1},\cdots,\varphi_{n})$.
To reach this goal, we first show that the quantity
$C_{2}(K,\varphi_{1},\cdots,\varphi_{n})$ satisfies an inequality
which is stated in the following claim.

\bcl\label{eq-sh-18} $C_{2}(K,\varphi_{1},\cdots,\varphi_{n})\leq
\big(C_{2}(K,\varphi_{1},\cdots,\varphi_{n})\big)^{1-K^{1/(1-n)}}\mu_1+\mu_2,$
 where $$\mu_1=\big(q(n,K)\big)^{1+K^{1/(1-n)}}
\int_{\mathbb{S}^{n-1}}|\eta-T(\theta_{\epsilon})|^{1-n+K^{2/(1-n)}}d\sigma(\eta),$$
$q(n,K)$ is from Theorem L, $\mu_2=\mu_3+\mu_4,$
\be\label{fg-1.1}\mu_3=K\frac{\|\varphi_{1}\|_{\infty}}{n}+K\sum_{k=2}^{m}\frac{\|\varphi_{k}\|_{\infty}}{n^{2}(n+2)}\left[\frac{n+4}{4n(n+2)}\right]^{k-2},\ee
 and
$$\mu_{4}=\left(\frac{n}{n+1}+\frac{1}{n}\right)\|\varphi_{1}\|_{\infty}+\sum_{k=2}^{m}
\left[\frac{\delta(n)}{2n}+\frac{1}{2n^{2}(n+2)}\right]\left[\frac{n+4}{4n(n+2)}\right]^{k-2}\|\varphi_{k}\|_{\infty}.$$
\ecl

To prove the claim, we need the following preparation. Firstly, we
prove that

\beq\label{sun-3}|D_{f}(T(\theta))|&\leq&K\bigg\{\int_{\mathbb{S}^{n-1}}\frac{|\varphi_{0}(T(\theta))-\varphi_{0}(\eta)|^{2}}{|\eta-T(\theta)|^{n}}d\sigma(\eta)
+\frac{\mu_3}{K}\bigg\} \eeq almost everywhere in $Q^{n-1}$, where
$\varphi_{0}=f|_{\mathbb{S}^{n-1}}$.

Since $f$ is $K$-quasiconformal mapping, by Claim \ref{Ca-1}, we see
that

\be\label{ck-1.1}
\lim_{r\rightarrow1^{-}}|D_{f}(S(r,\theta))|^{n}\leq\lim_{r\rightarrow1^{-}}
KJ_{f}(S(r,\theta))\ee exists almost everywhere in $Q^{n-1}$. By Lemma H, we obtain

\be\label{v-ck-1.2}\lim_{r\rightarrow1^{-}}\frac{\partial f\circ S}{\partial\theta_{1}}(r,\theta)\times\cdots\times\frac{\partial f\circ S}{\partial\theta_{n-1}}(r,\theta)
=\frac{\partial f\circ T}{\partial\theta_{1}}(\theta)\times\cdots\times\frac{\partial f\circ T}{\partial\theta_{n-1}}(\theta)\ee
exists almost everywhere in $Q^{n-1}$.
It follows from (\ref{v-ck-1.2}), Lemma J and
$$\frac{\partial f\circ S}{\partial\theta_{1}}(r,\theta)=rf'(S(r,\theta))\frac{\partial T}{\partial \theta_{j}}~(j\in\{1,\ldots,n-1\})$$ that

\be\label{v-ck-1.3}M_{x}(\theta)\leq\lim_{r\rightarrow1^{-}}|D_{f}(S(r,\theta))|^{n-1}M_{T}(\theta)=|D_{f}(T(\theta))|^{n-1}M_{T}(\theta),
\ee where $M_{x}(\theta)$ and $M_{T}(\theta)$ are defined in Lemma
\ref{lem-mck}. From (\ref{M-1.1}) in Lemma \ref{lem-mck},
(\ref{ck-1.1}) and (\ref{v-ck-1.3}), we infer that

\beqq|D_{f}(T(\theta))|^{n}&\leq& K J_{f}(T(\theta))\\&\leq&
K\frac{M_{x}(\theta)}{M_{T}(\theta)}\bigg\{\int_{\mathbb{S}^{n-1}}\frac{|\varphi_{0}(T(\theta))-\varphi_{0}(\eta)|^{2}}{|\eta-T(\theta)|^{n}}d\sigma(\eta)
+\frac{\mu_3}{K}\bigg\}\\
&\leq&K|D_{f}(T(\theta))|^{n-1}\bigg\{\int_{\mathbb{S}^{n-1}}\frac{|\varphi_{0}(T(\theta))-\varphi_{0}(\eta)|^{2}}{|\eta-T(\theta)|^{n}}d\sigma(\eta)
+\frac{\mu_3}{K}\bigg\} \eeqq almost everywhere in $Q^{n-1}$, which
yields that (\ref{sun-3}).

Secondly, we show that for any $\epsilon>0$, there exists
$\theta_{\epsilon}\in Q^{n-1}$ such that

\beq\label{vck-1.01}
C_{2}(K,\varphi_{1},\cdots,\varphi_{n})&\leq&(1+\epsilon)|D_{f}(T(\theta_{\epsilon}))|+\mu_{4}.\eeq
 Since
$$P[\varphi_{0}]=f-\sum_{k=1}^{m}(-1)^{k}G_{k}[\varphi_{k}]$$ is
harmonic, by Lemma I, we see that

\beqq
\left|D_{f}(x)-\sum_{k=1}^{m}(-1)^{k}D_{G_{k}[\varphi_{k}]}(x)\right|&=&|D_{P[\varphi_{0}]}(x)|\leq{\rm
ess ~sup}_{|t|=1}|D_{P[\varphi_{0}]}(t)|\\
&\leq&{\rm ess ~sup}_{|t|=1}|D_{f}(t)|+\sum_{k=1}^{m}{\rm ess
~sup}_{|t|=1}|D_{G_{k}[\varphi_{k}]}(t)|,
 \eeqq
which, together with Lemmas \ref{eq-jj-1} and \ref{eq-jj-2}, gives
that

\beq \label{vck-1.5}|D_{f}(x)|&\leq&{\rm ess
~sup}_{|t|=1}|D_{f}(t)|+\sum_{k=1}^{m}|D_{G_{k}[\varphi_{k}]}(x)|\\
\nonumber &&+\sum_{k=1}^{m}{\rm ess
~sup}_{|t|=1}|D_{G_{k}[\varphi_{k}]}(t)|\leq{\rm ess
~sup}_{|t|=1}|D_{f}(t)|+\mu_{4}.
 \eeq
Hence (\ref{vck-1.01}) follows from  (\ref{vck-1.5}) and Claim
\ref{Ca-1}.

For  $\theta\in Q^{n-1}$, let \be
\label{ttt-t}\Lambda(\theta)=\int_{\mathbb{S}^{n-1}}\frac{|\varphi_{0}(T(\theta))-\varphi_{0}(\eta)|^{2}}{|\eta-T(\theta)|^{n}}d\sigma(\eta).\ee
Then by Theorem L, we have

\beqq\Lambda(\theta_{\epsilon})&\leq&C_{2}(K,\varphi_{1},\cdots,\varphi_{n})^{1-K^{1/(1-n)}}\\
\nonumber
&&\times\int_{\mathbb{S}^{n-1}}|\eta-T(\theta_{\epsilon})|^{1-n+K^{2/(1-n)}}
\frac{|\varphi_{0}(T(\theta_{\epsilon}))-\varphi_{0}(\eta)|^{1+K^{1/(1-n)}}}{|\eta-T(\theta_{\epsilon})|^{K^{2/(1-n)}+K^{1/(1-n)}}}d\sigma(\eta)\\
\nonumber
&\leq&C_{2}(K,\varphi_{1},\cdots,\varphi_{n})^{1-K^{1/(1-n)}}\mu_1,
 \eeqq
which, together with (\ref{sun-3}) and (\ref{vck-1.01}), gives that

\beq\label{v-ck-1.8}
C_{2}(K,\varphi_{1},\cdots,\varphi_{n})&\leq&(1+\epsilon)|D_{f}(T(\theta_{\epsilon}))|+\mu_{4}\\
\nonumber&\leq&K(1+\epsilon)\Lambda(\theta_{\epsilon}) +\mu_3 (1+\epsilon)+ \mu_{4}\\
\nonumber
&\leq&K(1+\epsilon)C_{2}(K,\varphi_{1},\cdots,\varphi_{n})^{1-K^{1/(1-n)}}\mu_1\\
\nonumber && +\mu_3 (1+\epsilon)+ \mu_{4}. \eeq

By letting $\epsilon\rightarrow0^{+}$, we get from (\ref{v-ck-1.8})
that
$$C_{2}(K,\varphi_{1},\cdots,\varphi_{n})\leq KC_{2}(K,\varphi_{1},\cdots,\varphi_{n})^{1-K^{1/(1-n)}}\mu_1+\mu_2,$$
which yields that

$$C_{2}(K,\varphi_{1},\cdots,\varphi_{n})\leq(K\mu_1+\mu_2)^{K^{1/(n-1)}}.$$

\bcl\label{sun-t} If $\big(1-K^{1/(1-n)}\big)\mu_1<1$, then
$$C_{2}(K,\varphi_{1},\cdots,\varphi_{n})\leq \mu_5,$$ where
$\mu_5=\frac{K^{1/(1-n)}\mu_1+\mu_2}{1-\big(1-K^{1/(1-n)}\big)\mu_1}.$
\ecl
 The proof of this claim easily follows from \cite[Lemma 2.9]{K2}.

In the following, an upper bound of
$C_{2}(K,\varphi_{1},\cdots,\varphi_{n})$ will be established. By
Claims \ref{eq-sh-18} and \ref{sun-t}, we obtain that

\beqq\label{v-ck-1.10}C_{2}(K,\varphi_{1},\cdots,\varphi_{n})\leq
C_{3},\eeqq where

\be\label{v-ck-1.11} C_{3}=\begin{cases}
\displaystyle (K\mu_1+\mu_2)^{K^{1/(n-1)}}, &\mbox{ if }\, \big(1-K^{1/(1-n)}\big)\mu_1\geq1,\\
\displaystyle \min\left\{\mu_5,
(K\mu_1+\mu_2)^{K^{1/(n-1)}}\right\}, &\mbox{ if }\,
\big(1-K^{1/(1-n)}\big)\mu_1<1.
\end{cases}\ee

In the following, we will  break $C_{3}$ down into the form we need.
By (\ref{v-ck-1.11}), we have

$$C_3=\begin{cases}
\displaystyle M^{\ast}_{1}, &\mbox{ if }\, \big(1-K^{1/(1-n)}\big)\mu_1\geq 1,\\
\displaystyle \min\big\{M^{\ast}_{1}, M^{\ast}_{2}\big\}, &\mbox{ if
}\, \big(1-K^{1/(1-n)}\big)\mu_1< 1,
\end{cases}$$
where
$M^{\ast}_{1}=M_{1}^{'}(n,K)+N_{1}^{'}(K,\varphi_{1},\cdots,\varphi_{n})$,
$M^{\ast}_{2}=M_{1}^{''}(n,K)+N_{1}^{''}(K,\varphi_{1},\cdots,\varphi_{n})$,
$M_{1}^{'}(n,K)=(K\mu_1)^{K^{1/(n-1)}}$,
$M_{1}^{''}(n,K)=\frac{K^{1/(1-n)}\mu_1}{1-\big(1-K^{1/(1-n)}\big)\mu_1}$,
$N_{1}^{'}(K,\varphi_{1},\cdots,\varphi_{n})=(K\mu_1+\mu_2)^{K^{1/(n-1)}}-(K\mu_1)^{K^{1/(n-1)}}$,
and
$N_{1}^{''}(K,\varphi_{1},\cdots,\varphi_{n})=\frac{\mu_2}{1-\big(1-K^{1/(1-n)}\big)\mu_1}.$

Let
$$  M_{1}(n,K)=\begin{cases}
\displaystyle M_{1}^{'}(n,K), &\mbox{ if }\, \big(1-K^{1/(1-n)}\big)\mu_1\geq1,\\
\displaystyle M_{1}^{''}(n,K), &\mbox{ if }\,
\big(1-K^{1/(1-n)}\big)\mu_1<1 ~\mbox{and}~M_{1}^{\ast}\geq
M_{2}^{\ast},\\
\displaystyle M_{1}^{'}(n,K), &\mbox{ if }\,
\big(1-K^{1/(1-n)}\big)\mu_1<1 ~\mbox{and}~M_{1}^{\ast}\leq
M_{2}^{\ast}
\end{cases}$$
and
$$ N_{1}(K,\varphi_{1},\cdots,\varphi_{n})=\begin{cases}
\displaystyle N_{1}^{'}(K,\varphi_{1},\cdots,\varphi_{n}), &\mbox{ if }\, \big(1-K^{1/(1-n)}\big)\mu_1\geq1,\\
\displaystyle N_{1}^{''}(K,\varphi_{1},\cdots,\varphi_{n}), &\mbox{
if }\, \big(1-K^{1/(1-n)}\big)\mu_1<1 ~\mbox{and}~M_{1}^{\ast}\geq
M_{2}^{\ast},\\
\displaystyle N_{1}^{'}(K,\varphi_{1},\cdots,\varphi_{n}), &\mbox{
if }\, \big(1-K^{1/(1-n)}\big)\mu_1<1 ~\mbox{and}~M_{1}^{\ast}\leq
M_{2}^{\ast}.
\end{cases}$$


It follows from the facts
$$\lim_{K\rightarrow1^{+}}M_{1}(n,K)=1\;\;\mbox{and}\;\; \lim_{\|\varphi_{1}\|_{\infty}\rightarrow0^{+},\cdots,\|\varphi_{n}\|_{\infty}\rightarrow0^{+}}N_{1}(K,\varphi_{1},\cdots,\varphi_{n})=0$$
that these two constants are what we need. The proof of this theorem
is complete. \qed


\bigskip
\section{Acknowledgments}
This research  was partly supported by the National Science
Foundation of China (grant no. 12071116), the Hunan Provincial Natural Science Foundation of China (No. 2022JJ10001),
 the Key Projects of Hunan Provincial Department of Education (grant no. 21A0429);
 the Double First-Class University Project of Hunan Province
(Xiangjiaotong [2018]469),  the Science and Technology Plan Project of Hunan
Province (2016TP1020), and the Discipline Special Research Projects of Hengyang Normal University (XKZX21002)



\end{document}